\RequirePackage{fix-cm}
\documentclass[numbook,smallcondensed]{svjour3}     % onecolumn (ditto)
\smartqed  % flush right qed marks, e.g. at end of proof
\journalname{}
\usepackage{amsmath}
\usepackage{amssymb}
\usepackage{amsfonts}
\usepackage{graphicx}
\usepackage{algorithm}
\usepackage{algorithmic}
\usepackage{cite}
\usepackage{enumitem}
\usepackage{bm}
\usepackage[caption=false]{subfig}
\usepackage{tikz}
\usepackage{pgfplots}
\usepackage{multirow}
\usepackage{lscape}
\usepackage{textcomp}
\usepackage{chngcntr}
\pgfplotsset{compat=newest}
\usetikzlibrary{plotmarks}

\makeatletter
\let\cl@chapter\undefined
\makeatletter

\counterwithout{figure}{section}
\counterwithout{table}{section}
\numberwithin{algorithm}{section}

\usepackage[%
  colorlinks = true,
  linkcolor  = linkBlue,
  citecolor  = linkRed,
  urlcolor   = linkPurple
]{hyperref}
\usepackage{doi}
\usepackage[nameinlink,capitalise]{cleveref}
\crefformat{equation}{(#2#1#3)}

\newcommand{\bx}[1]{\ensuremath{\mathbf{#1}}}
\newcommand\R{\ensuremath{\mathbb{R}}}
\newlength\fheight
\newlength\fwidth
\definecolor{linkBlue}{HTML}{0055C9}
\definecolor{linkRed}{HTML}{FF1A24}
\definecolor{linkPurple}{HTML}{6200D9}

\definecolor{matlabBlue}{HTML}{0072BD}
\definecolor{matlabGreen}{HTML}{77AC30}
\definecolor{matlabBrown}{HTML}{A2142F}

\newcommand{\RomanNumeralCaps}[1]
{\MakeUppercase{\romannumeral #1}}

\begin{document}

\title{A Training Set Subsampling Strategy for the Reduced Basis Method\thanks{The first author is supported by the International Max Planck Research School for Advanced Methods in Process and Systems Engineering (IMPRS ProEng).}}

\titlerunning{A Training Set Subsampling Strategy for the Reduced Basis Method}        % if too long for running head

\author{Sridhar Chellappa \and Lihong Feng \and \newline Peter Benner
}


\institute{S. Chellappa \at
              Max Planck Institute for Dynamics of Complex Technical Systems, Sandtorstra\ss e 1, 39106 Magdeburg, Germany \\
              \email{chellappa@mpi-magdeburg.mpg.de}           %  \\
           \and
           L. Feng \at
              Max Planck Institute for Dynamics of Complex Technical Systems, Sandtorstra\ss e 1, 39106 Magdeburg, Germany \\
							\email{feng@mpi-magdeburg.mpg.de}
					 \and
					 P. Benner \at
							Max Planck Institute for Dynamics of Complex Technical Systems, Sandtorstra\ss e 1, 39106 Magdeburg, Germany \\
							\email{benner@mpi-magdeburg.mpg.de}
}

\date{Received: date / Accepted: date}

\maketitle

\begin{abstract}
	We present a subsampling strategy for the offline stage of the Reduced Basis Method. The approach is aimed at bringing down the considerable offline costs associated with using a finely-sampled training set. The proposed algorithm exploits the potential of the pivoted QR decomposition and the discrete empirical interpolation method to identify important parameter samples. It consists of two stages. In the first stage, we construct a low-fidelity approximation to the solution manifold over a fine training set. Then, for the available low-fidelity snapshots of the output variable, we apply the pivoted QR decomposition or the discrete empirical interpolation method to identify a set of sparse sampling locations in the parameter domain. These points reveal the structure of the parametric dependence of the output variable. The second stage proceeds with a subsampled training set containing a by far smaller number of parameters than the initial training set. Different subsampling strategies inspired from recent variants of the empirical interpolation method are also considered. Tests on benchmark examples justify the new approach and show its potential to substantially speed up the offline stage of the Reduced Basis Method, while generating reliable reduced-order models.
\end{abstract}

\section{Introduction}
\label{sec:intro}
Mathematical models in many areas of science and engineering are described by parametric ordinary differential equations (ODEs) or partial differential equations (PDEs). The analytical solution of such ODEs and PDEs is usually intractable, therefore engineers prefer detailed numerical simulations of the models to obtain useful insights into the underlying processes. However, a detailed simulation of the model necessitates a very fine discretization of the governing ODEs or PDEs, both in space and time. The resulting discretized system of equations is referred to as the high-fidelity model or the full-order model (FOM). Numerical solutions of the FOM are computationally expensive to obtain, especially when computing resources are limited or when repeated simulations of the FOM are needed. The latter scenario is common in a multi-query task, such as uncertainty quantification, optimization, etc. With the goal of speeding up simulations, the area of model order reduction (MOR) has gained popularity in the last decades \cite{morAnt05,AntBG20,morBenetal21,morBenCOetal17,morBenGW15}. Different MOR algorithms have been proposed to obtain surrogate models for the FOMs. These surrogates are often referred to as reduced-order models (ROMs) since they have a reduced number of degrees of freedom compared to the FOMs. The simulation results of a suitably developed ROM are indistinguishable from those of the FOM, yet, the ROM incurs only a fraction of the computational cost spent on the FOM. This is a big advantage in real-time simulation or multi-query scenarios. Generally, there exist two classes of ROMs:  projection-based ROMs \cite{morAnt05} which are obtained by projecting the original model onto an approximation subspace, and data-driven ROMs, which are obtained through a data fit \cite{morPehW16}. In the following, we shall limit ourselves to projection-based ROMs. 

The Reduced Basis Method (RBM) is a popular, projection-based technique to obtain ROMs of parametric PDEs. The ansatz underlying the RBM is that the solution manifold of the parametric FOM can be well-approximated by a low-dimensional subspace. To achieve efficient online simulation, considerable efforts are made to explore the parameter space and collect solution snapshots in the offline stage. The greedy algorithm explores the parameter space based on suitable \textit{a posteriori} error estimators. At each iteration of the greedy algorithm, the parameter maximizing the estimated error over a training set is selected and the corresponding FOM is solved to obtain new snapshots. The training set is, in essence, a discrete representation of the parameter domain. The snapshots are used to construct a basis for the low-dimensional approximation space, followed by a (Petrov-)Galerkin projection to obtain the ROM. We refer to the recent books \cite{morHesRS16, morQuaMN16} and the survey works \cite{morCheQR17,morQuaRM11,morRozHP08} for a detailed exposition of the theoretical framework underpinning the RBM and also for applications of the RBM in a variety of problems.

The choice of the training set is critical for the success of the RBM. A poorly-sampled training set can result in an inadequate representation of all the solution modes, causing the ROM to fail to meet the desired tolerance criterion for a parameter not present in the training set. Therefore, it is common practice to adopt a finely-sampled training set. However, the computational cost of the offline stage scales with the cardinality of the training set, which becomes high for problems with high-dimensional parameter space. Therefore, a more efficient sampling strategy is desired. 

Many works have attempted to address the issue of optimal training set sampling. Notable among them are: the Multi-Stage Greedy algorithm from \cite{morSen08} and the Adaptively Enriching Greedy algorithm in \cite{morHesSZ14}. In the former, the author suggests performing a set of greedy algorithms over randomly sampled training sets; then the resulting ROM is tested over a much larger random training set and the greedy algorithm is re-run on those points failing to meet the tolerance criterion in the larger training set. In the latter work, the authors propose a saturation criterion which is used to systematically remove parameters from a randomly-sampled training set. New random parameters are then added to the current training set. A larger training set serves as a safety check mechanism at every iteration. However, it may not be efficient, in general, to estimate a robust saturation criterion. The authors of \cite{morEftKP11,morHaaDO11} propose a localized RBM approach, where a hierarchical tree-based partitioning of the parameter domain is done and separate ROMs for each partition are generated. In \cite{morHesZ16}, the authors consider a two-stage approach that uses the ANOVA expansion together with parameter domain decomposition to address training set complexity. The work \cite{morMadS13} considers an anisotropic sampling of the parameter domain using an empirical norm derived from the truncated Hessian of the solution vector with respect to the parameter. No explicit partition of the parameter domain is considered. However, the basis vectors are determined in the online stage. Moreover, the method needs to compute the Hessian at each point in the training set in order to define a distance metric which is subsequently used to add more samples to the training set. The calculation of the Hessian can be very expensive, especially for non-stationary problems. The method proposed in \cite{morJiaCN17} makes only a subset of the finely sampled training set active at a given iteration of the greedy algorithm. A recent extension of this work \cite{morJiaC20} proposes hybrid strategies combining the ideas from \cite{morHesSZ14,morSen08}. Different strategies are proposed to identify the set of active parameters. The works \cite{morCheFB19,morTaineA15} propose a cheap surrogate model for a certain error estimator, based on Kriging and radial basis functions, respectively, and use it to obtain the estimated error for any parameter in a fine training set. A sparse grid-based construction for the training set is suggested in \cite{morPehZB13}.

The authors of \cite{morCheG19} perform an eigendecomposition of the Hessian matrix of the output variable with respect to the parameter to identify a small subspace of the high-dimensional parameter space by truncation. The parameters that constitute the training set for the RBM are then sampled from the identified eigenspace. In \cite{morTezBR18}, the so-called \emph{active subspace} \cite{morConDW14} of the parameter space is identified by relying on gradient information of the output with respect to the parameter. Both these works are limited to scalar-valued outputs and steady-state problems.

Most of the existing work related to adaptive training set sampling focuses on steady-state or quasi steady-state problems. To the best of our knowledge, only \cite{morCheFB19,morEftKP11,morHaaDO11,morTaineA15} address training set adaptivity for time-dependent problems. The works \cite{morEftKP11,morHaaDO11} propose a localization strategy that involves constructing multiple ROMs over local parameter domains, while \cite{morCheFB19,morTaineA15} consider adaptively enriching a coarse training set by observing a cheap error surrogate over a fine training set.

In this work, we present a goal-oriented training set sampling strategy that relies on the output quantity of interest (QoI). We aim at identifying the structure of the parameter dependency of the output through the empirical interpolation algorithm \cite{morBarMNetal04,morChaS10} or a pivoted QR decomposition and utilize this information to find out the parameter importance. Our proposed method is applicable to both steady-state and time-dependent problems with vector-valued outputs. Our central contribution is a two-stage algorithm to control the cardinality of the training set. In the first stage, a coarse RB approximation of the problem is obtained using a fine training set. Then, an approximate output snapshot matrix is derived by time integrating the coarse ROM at all the parameter samples in the fine training set. We apply the pivoted QR decomposition or, alternatively, the discrete empirical interpolation method (\textsf{DEIM}) (or its variants) to the approximate output snapshot matrix. This procedure identifies regions of the parameter space that have a greater contribution to the current RB approximation space. In the second stage, the fine training set is subsampled based on the parameter distribution identified using the pivoted QR decomposition or the \textsf{DEIM} algorithm and leads to a subsampled coarse training set. The RBM is continued over the coarse training set, until a targeted error tolerance is met.

The paper is organized as follows. In \Cref{sec:probsetting} we describe the preliminaries including the problem setting for the proposed methodology, the RBM and the related hyper-reduction techniques. In \Cref{sec:trngset} we detail the issue of training set sampling for the RBM and present our main algorithm for efficient training set subsampling. \Cref{sec:numerics} is dedicated to numerical examples through which we illustrate various aspects of the proposed subsampling strategy and demonstrate the speedup it offers for two numerical benchmark problems. We conclude by summarizing the proposed method and highlighting possible research directions for the near future. 
Throughout this work, we have used $\operatorname{MATLAB}\textsuperscript{\textregistered}$ notation in the presentation of algorithms and numerical experiments.

\section{Preliminaries}
\label{sec:probsetting}
In this section, we present the continuous problem and the discretized  system that the proposed subsampling strategy is valid for. Then, we briefly review the RBM and the associated issue of training set sampling. Afterwards, some hyper-reduction algorithms are reviewed in order to introduce our proposed subsampling algorithms.
\subsection{High-fidelity Models}
A wide variety of physical and engineering phenomena are modelled via PDEs. Consider the spatial domain $\Omega \subset \R^{d}$ with $(d = 1, 2, 3)$. Let a model of PDEs defined in $\Omega$ be denoted by
  \begin{equation}\label{eq:pde}
    \mathcal{L}\left(\bx{v}, \bx{w}, t, \bm{\mu} \right) = \bm{0},
  \end{equation}
where $\bx{v}$ is the (vector-valued) state variable and describes the particular physical quantity the PDE models,  $\bx{w} \in \Omega$ is the spatial variable, $0 \leq t \leq T$ denotes the time and $\bm{\mu} \in \mathcal{P} \subset \R^{p}$ defines the parameters. The above form is a general description of any time-dependent or steady-state problem with or without parameter variations. The output of the model is usually a function of the solution $\bx{v}$ and the parameter $\bm{\mu}$. After numerical discretization in space and time, we get
  \begin{equation}\label{eq:odess}
  	\begin{aligned}
	    \bx{E}\, \bx{x}(t^{k},\bm{\mu}) &= \bx{A}\, \bx{x}(t^{k-1},\bm{\mu}) + \bx{f}\left (\bx{x}(t^{k-1},\bm{\mu}), \bm{\mu}\right) + \bx{B}\, \bx{u}(t^{k-1},\bm{\mu}),\\
	    \bx{y}\left(\bx{x}(t^{k},\bm{\mu}), \bm{\mu}\right) &= \bx{C}\, \bx{x}(t^{k},\bm{\mu}).
  	\end{aligned}
  \end{equation}
Here, $\bx{x}(t^{k}, \bm{\mu}) \in \R^{n}$ is the state vector, $\bx{u}(t^{k}, \bm{\mu}) \in \R^{m}$ is the input vector and $\bx{y}\left(\bx{x}(t^{k},\bm{\mu}), \bm{\mu}\right) \in \R^{q}$ is the output or quantity of interest. Further, $\bx{E}, \bx{A} \in \R^{n \times n}$, $\bx{B} \in \R^{n \times m}$ is the input matrix, $\bx{C} \in \R^{q \times n}$ is the output matrix and $\bx{f}\left(\bx{x}(t^{k},\bm{\mu}), \bm{\mu}\right) \in \R^{n}$ models the nonlinearity associated with the system. 
\begin{remark}
The system matrices ($\bx{E}, \bx{A}, \bx{B} \, \text{and}\, \bx{C}$) can also be time- and/or parameter-dependent. However, we have not made this dependence explicit for the sake of keeping the notations concise. For the case of steady-state problems, the time dependence of the state, input and output vectors and system matrices vanishes and the system simply reads
	\begin{equation}\label{eq:odess_static}
		\begin{aligned}
		\bx{E}\, \bx{x}(\bm{\mu}) &= \bx{f}\left (\bx{x}(\bm{\mu}), \bm{\mu}\right) + \bx{B}\, \bx{u}(\bm{\mu}),\\
		\bx{y}(\bx{x}(\bm{\mu}), \bm{\mu}) &= \bx{C}\, \bx{x}(\bm{\mu}).
		\end{aligned}
	\end{equation}	
\end{remark}
The discretized system in \cref{eq:odess} or \cref{eq:odess_static} is also called the FOM and has a large number of degrees of freedom, i.e., $n$ is very large. The proposed subsampling strategy is applicable to both \cref{eq:odess,eq:odess_static}.
\subsection{Reduced Basis Method and the Training Set}
\label{sec:rbm}
The Reduced Basis Method relies on the observation that the solution manifold can be well-approximated by a \emph{low-dimensional} subspace $\mathcal{V}$. Let $\left[ \bx{v}_{1}, \, \bx{v}_{2}, \, \ldots, \, \bx{v}_{r} \right] =: \bx{V} \in \R^{n \times r}$ be a basis of the subspace $\mathcal{V}$. The approximated solution for any parameter is obtained by considering the ansatz $\bx{x}(t^{k}, \bm{\mu}) \approx \widetilde{\bm{x}}(t^{k}, \bm{\mu}) = \sum_{i = 1}^{r} z_{i} \bx{v}_{i}$. The parameter-dependent coefficients $\bx{z} = [z_{1}, \, z_{2}, \, \ldots, \, z_{r}]^{\textsf{T}}$ are obtained by solving the ROM
  \begin{equation}\label{eq:odess_r}
  	\begin{aligned}
	    \bx{E}_{r}\, \bx{z}(t^{k},\bm{\mu}) &= \bx{A}_{r}\, \bx{z}(t^{k-1},\bm{\mu}) + \bx{f}_{r}\left(\bx{V} \bx{z}(t^{k-1},\bm{\mu}), \bm{\mu}\right) + \bx{B}_{r}\, \bx{u}(t^{k-1},\bm{\mu}),\\
	    \widetilde{\bx{y}}\left(\bx{z}(t^{k},\bm{\mu}), \bm{\mu}\right) &= \bx{C}_{r}\, \bx{z}(t^{k},\bm{\mu}),
  	\end{aligned}
  \end{equation}
derived through Galerkin projection of the FOM onto the low-dimensional subspace $\mathcal{V}$.
The reduced system matrices $\bx{E}_{r},\, \bx{A}_{r} \in \R^{r \times r}$, $\bx{B}_{r} \in \R^{r \times m}$ and $\bx{C}_{r} \in \R^{q \times r}$ are obtained as $\bx{E}_{r} := \bx{V}^{\textsf{T} } \bx{E} \bx{V}$, $\bx{A}_{r} := \bx{V}^{\textsf{T} } \bx{A} \bx{V}$, $\bx{B}_{r} := \bx{V}^{\textsf{T} } \bx{B}$ and $\bx{C}_{r} := \bx{C} \bx{V}$, respectively. Finally, $\bx{f}_{r} := \bx{V}^{\textsf{T}} \bx{f}\left( \bx{V} \bx{z}(t^{k-1},\bm{\mu}), \bm{\mu} \right)$.

The greedy algorithm or the POD-greedy (Proper Orthogonal Decomposition-greedy) algorithm are the most popular techniques for constructing the RBM approximation space for the steady-state system \cref{eq:odess_static} and the time-dependent system \cref{eq:odess}, respectively. In order to initialize the greedy algorithm, a training set $\Xi_{\textnormal{train}}$ is given \emph{a priori}, from which parameter samples are iteratively selected, so that the corresponding solution vector is iteratively added to the basis matrix $\bx{V}$. We summarize the (POD-)greedy algorithm in \cref{alg:rb}, which takes both the steady-state case and the time-dependent case into consideration. In Step 8 of \cref{alg:rb}, the snapshot matrix $\bm{\mathfrak{X}}$ for the time-dependent case consists of the solution vector at discretized time instances $\{t^{i}\}_{i=0}^{K}$ given by:
	\begin{displaymath}
		\bm{\mathfrak{X}}(\bm{\mu}) = \left[ \bx{x}(t^{0}, \bm{\mu}) \cdots \bx{x}(t^{K}, \bm{\mu}) \right ] \in \R^{n \times N_{t}},
	\end{displaymath}
where $N_{t} := (K + 1)$. $\bx{U}_{\bar{\bm{\mathfrak{X}}}}$ in Step 10 is the left singular vector matrix obtained from the singular value decomposition (SVD) of $\bar{\bm{\mathfrak{X}}}$, i.e., $\bar{\bm{\mathfrak{X}}} = \bx{U}_{\bar{\bm{\mathfrak{X}}}} \Sigma_{\bar{\bm{\mathfrak{X}}}} \bx{V}_{\bar{\bm{\mathfrak{X}}}}^{\textsf{T}}$. Here, $\Sigma_{\bar{\bm{\mathfrak{X}}}}$ contains all the non-zero singular values of $\bar{\bm{\mathfrak{X}}}$: $\sigma_{1} \geq \sigma_{2} \geq \cdots \geq \sigma_{r_{X}} \geq 0$. If no alternative decision criterion is used, $r_{\text{POD}}$ is usually taken to be $1$. In Step 13, $\Delta(\bm{\mu})$ denotes an error estimator for the error in approximating the state variable or the output variable. The error estimator should be much cheaper to compute than the true error. For the sake of clarity, the sketched algorithm is the basic version of the RBM. Several enhancements in the form of primal-dual error estimation, hyper-reduction, adaptive basis construction, etc. exist \cite{morBenEEetal18, morCheFB20,morGreMetal07,morGreP05,morZhaetal15}. 
\renewcommand{\algorithmicrequire}{\textbf{Input:}}
\renewcommand{\algorithmicensure}{\textbf{Output:}}
\renewcommand{\algorithmiccomment}[1]{// #1}
\begin{algorithm}[t!]
\caption{Reduced Basis Method (RBM)}
\label{alg:rb}
\begin{algorithmic}[1]
\REQUIRE{Training set ($\Xi$), tolerance for the ROM ($\texttt{tol}$), maximum iterations ($n_{\text{max}}$).}
\ENSURE{$\bx{V}$.}
\STATE{Initialize $\bx{V} = [\,\,]$, $\epsilon = 1 + \texttt{tol}$, \texttt{iter} = 1.}
\STATE Select $\bm{\mu}^{*}_{(1)}$ randomly from $\Xi_{\text{train}}$. 
\WHILE{$\epsilon > \texttt{tol}\,\,\& \& \,\, \texttt{iter} \leq n_{\text{max}}$}
	\IF{\textbf{steady-state}}
		\STATE{Solve FOM \cref{eq:odess_static} at $\bm{\mu}^{*}_{(\texttt{iter})}$ and obtain solution $\bx{x}(\bm{\mu}^{*}_{(\texttt{iter})})$.}
		\STATE{Set $\bx{V} \leftarrow{} \texttt{orth}\bigg(\left[ \bx{V} \,\, \bx{x}(\bm{\mu}^{*}_{(\texttt{iter})})\right]\bigg)$. \hfill \% \textit{\texttt{orth} means orthogonalizing $\bx{x}(\bm{\mu}^{*}_{(\texttt{iter})})$ against the column vectors in the current matrix $\bx{V}$.
		}}		
	\ELSE{}
		\STATE{Solve FOM \cref{eq:odess} at $\bm{\mu}^{*}_{(\texttt{iter})}$ and obtain snapshot matrix $\bm{\mathfrak{X}}$.}
		\STATE Set $\bar{\bm{\mathfrak{X}}} = \bm{\mathfrak{X}} - \bx{V} \bx{V}^{\textsf{T}} \bm{\mathfrak{X}}$.
		\STATE{Enrich $\bx{V}$ with $r_{\text{POD}}$ left singular vectors of  $\bar{\bm{\mathfrak{X}}}$ as
		
		$\bx{V} \leftarrow{} \texttt{orth}\left(\left[ \bx{V} \,\, \bx{U}_{\bar{\bm{\mathfrak{X}}}}(:\,,\,1 \, : \, r_{\text{POD}}) \right]\right)$.}
	\ENDIF		
	\STATE{\texttt{iter} = \texttt{iter} + 1.}
	\STATE{Set $\bm{\mu}^{*}_{(\texttt{iter})} = \operatorname{arg} \operatorname{max} \limits_{\bm{\mu} \in \Xi_{\text{train}}} \Delta(\bm{\mu})$.}
	
	\STATE{$\epsilon = \Delta(\bm{\mu}^{*}_{(\texttt{iter})})$.}
	
\ENDWHILE
\end{algorithmic}
\end{algorithm}

It is noticed that the standard greedy algorithm (\cref{alg:rb}) does not address the issue of properly choosing $\Xi_{\textnormal{train}}$. When the parameter space dimension is high ($p \gg 2$) or when the number of time steps $N_{t}$ is large, the overall computational cost for the greedy algorithm to construct $\bx{V}$ can be substantial. The technique considered in this work aims at reducing the offline cost by subsampling a fine training set. We detail this in \Cref{sec:trngset}. Before that, since it will be used later, we briefly review two hyper-reduction procedures --- the \textsf{DEIM} algorithm (and its variants) and the \textsf{Gappy-POD} method. 
\subsection{Discrete Empirical Interpolation Method}
	\label{subsec:deim}
\begin{algorithm}[t!]
\caption{Discrete Empirical Interpolation Method (\textsf{DEIM}) as in \cite{morChaS10}}
\label{alg:deim}
\begin{algorithmic}[1]
\REQUIRE{Snapshots of the nonlinear vector ($\bm{\mathfrak{F}}$), tolerance for the SVD ($\epsilon_{\text{SVD}}$).}
\ENSURE{$\bx{U}$, $\bx{P}$ and $\bx{I}$.}
\STATE{Perform SVD of $\bx{\mathfrak{F}}$; collect the first $\ell$ left singular vectors $\{ \bx{u}_{j} \}_{j = 1}^{\ell} \subset \R^{n}$.}
\STATE {Set $[\sim \, , p_{1}]\, =\, \operatorname{max}(|\bx{u}_{1}|)$. \hfill \% \textit{$|\cdot|$ denotes the absolute value.}}
\STATE Set $\bx{U} = \bx{u}_{1}$, $\bx{I} = p_{1}$ and $\bx{P} = [e_{p_{1}}]$.
\FOR{$i = 2\,\, \text{to}\,\, \ell$}
	\STATE{Solve $\bx{P}^{\textsf{T}} \bx{U} \bm{\alpha} = \bx{P}^{\textsf{T}} \bx{u}_{i}$}.
	\STATE{Define residual $\bx{r} = \bx{u}_{i} - \bx{U} \bm{\alpha}$}.
	\STATE{Set $[\sim \, , p_{i}]\, =\, \operatorname{max}(|\bx{r}|)$}.
	\STATE{Update $\bx{U} := [ \bx{U}\,\, \bx{u}_{i} ]$, $\bx{I} := [\bx{I}\,\, p_{i}]$ and $\bx{P} := [ \bx{P}\,\, e_{p_{i}} ]$.}
\ENDFOR
\end{algorithmic}
\end{algorithm}
The \textsf{DEIM} algorithm \cite{morBarMNetal04,morChaS10,morGreMetal07} was introduced in the context of MOR, for efficient calculation of nonlinear or nonaffine terms of the ROM. The algorithm proceeds by collecting snapshots of the nonlinear vector in the FOM in \cref{eq:odess} or \cref{eq:odess_static} for different values of $\bm{\mu} \in \Xi_{\text{train}}$ given by
	\begin{displaymath}
		\bm{\mathfrak{F}} = \left[ \bx{f}\left (\bx{x}, \bm{\mu}_{1}\right) \cdots \bx{f}\left (\bx{x}, \bm{\mu}_{N_{\text{train}}}\right) \right ] \in \R^{n \times N_{\text{train}}}.
	\end{displaymath}
In order to avoid any $n$-dependent operations to evaluate the nonlinear vector involved in simulating the ROM, \textsf{DEIM} considers an approximation of the nonlinear vector given by
	\begin{equation}\label{eq:deim}
		\bx{f}\left( \bx{x}, \bm{\mu} \right) \approx \widetilde{\bx{f}}\left( \bx{x}, \bm{\mu} \right) = \bx{U} \bm{\alpha},
	\end{equation}
where $\bm{\alpha} \in \R^{\ell}$. The columns of $\bx{U} \in \R^{n \times \ell}$ are the interpolation basis vectors obtained via SVD of $\bm{\mathfrak{F}}$. The number of basis vectors $\ell$ can be determined through an information-theoretic criterion that depends on the relative energy content of the singular values $\{\sigma_{i}\}_{i=1}^{r_{X}}$ and reads $\frac{\sum\limits_{i=\ell+1}^{r_{X}} \sigma_i}{\sum\limits_{i=1}^{r_{X}} \sigma_i}~<~\epsilon_{\text{SVD}}$ where $\epsilon_{\text{SVD}}$ is a user-defined tolerance. Since \cref{eq:deim} is an overdetermined system with $n \gg \ell$, \textsf{DEIM} solves for $\bm{\alpha}$ by selecting $\ell$ rows from $\bx{U}$ and enforces interpolation as below:
	\begin{equation}\label{eq:deim_interp}
		\bx{P}^{\textsf{T}} \bx{f}\left( \bx{x}, \bm{\mu} \right) = \bx{P}^{\textsf{T}} \bx{U} \bm{\alpha}.
	\end{equation}
The indices of the rows where interpolation is enforced are given by $\{p_{1}, \ldots, p_{\ell}\}$. The $i^{\text{th}}$ column denoted as $e_{p_{i}} \in \R^{n}$ of the matrix $\bx{P} \in \R^{n \times \ell}$ is essentially the $i^{\text{th}}$ canonical unit vector with zeros at all but the $p_{i}^{\text{th}}$ entry. A greedy procedure, shown in \cref{alg:deim}, is used to identify $\bx{P}$. The algorithm ensures that $\left( \bx{P}^{\textsf{T}} \bx{U} \right)$ is nonsingular, so that $\widetilde{\bx{f}}\left( \bx{x}, \bm{\mu} \right)$ at any parameter $\bm{\mu}$ is given by
	\begin{displaymath}
		\widetilde{\bx{f}}\left( \bx{x}, \bm{\mu} \right) = \bx{U} \left( \bx{P}^{\textsf{T}} \bx{U} \right)^{-1} \bx{P}^{\textsf{T}} \bx{f}\left( \bx{x}, \bm{\mu} \right),
	\end{displaymath}
where $\bx{U} \left( \bx{P}^{\textsf{T}} \bx{U} \right)^{-1} \in \R^{n \times \ell}$ can be precomputed and stored. Moreover, the original nonlinear vector needs to be evaluated only at $\ell$ points, limiting the cost of evaluating the nonlinear vector to order $\ell$ and independent of $n$. The error in approximating the nonlinear vector is quantified as
	\begin{equation}
		\left \lVert \bx{f}(\bx{x}, \bm{\mu}) - \widetilde{\bx{f}}(\bx{x}, \bm{\mu}) \right \rVert_{2} \leq \lVert \left(\bx{P}^{\textsf{T}} \bx{U}\right)^{-1}  \rVert_{2} \lVert \bx{f}(\bx{x}, \bm{\mu}) - \bx{U} \bx{U}^{\textsf{T}} \bx{f}(\bx{x}, \bm{\mu}) \rVert_{2}.
	\end{equation}
The greedy choice of the interpolation indices in the \textsf{DEIM} algorithm is geared towards minimizing the term $\lVert \left(\bx{P}^{\textsf{T}} \bx{U}\right)^{-1}  \rVert_{2}$ appearing in the error bound. At each iteration, the new sampling point is chosen as the one resulting in the maximum reduction in $\lVert \left(\bx{P}^{\textsf{T}} \bx{U}\right)^{-1}  \rVert_{2}$.

Several variants of the \textsf{DEIM} algorithm have been proposed \cite{morDrmG16,morNegMA15,morPetal14,morPehDG18,morPehW15}. Among those, the \textsf{QDEIM} algorithm from \cite{morDrmG16} and the \textsf{DEIM}-based oversampling strategies proposed in \cite{morPehDG18} shall be of particular interest to our discussion in \Cref{sec:trngset}.
\subsubsection{QDEIM}
	\label{subsec:qdeim}
	In contrast to the original \textsf{DEIM} algorithm, the \textsf{QDEIM} approach from \cite{morDrmG16} relies on a column-pivoted QR decomposition to identify the interpolation points. This is different from the sequential, greedy choice of interpolation points in \textsf{DEIM} (see Steps 4 - 9 in \cref{alg:deim}).
	\textsf{QDEIM} is proven to have a sharper error bound and is also computationally more efficient and straightforward to implement. 
\subsubsection{KDEIM}
The \textsf{K} in  the \textsf{KDEIM} algorithm refers to the \texttt{k-means} clustering algorithm. The \texttt{k-means} algorithm is applied to the matrix $\bx{U}$ of (truncated) left singular vectors obtained from SVD of the snapshot matrix $\bm{\mathfrak{F}}$, then rows with similar response are assigned to the same cluster. The standard \texttt{k-means} objective function is recast as a relaxed trace maximization problem which is then solved using the QR decomposition. We refer the interested reader to \cite{morPehDG18} for a deeper discussion. Another early work to consider QR decomposition based clustering in the context of MOR was \cite{morMliGB15}, which used it for reducing networked multi-agent systems.
\subsubsection{Gappy-POD}
The number of interpolation points of the \textsf{DEIM} algorithm and its variants discussed so far equals to the number $\ell$ of interpolation basis vectors. However, in many cases it is beneficial to consider $m > \ell$ interpolation points. The \textsf{Gappy-POD} method and other related approaches fall into this category \cite{morCaretal11,Eve95}. The coefficient matrix $\bx{P}^{\textsf{T}}\bx{U}$ of the linear system in \cref{eq:deim_interp} is no longer square and, therefore, does not possess a unique inverse; instead it is solved using the pseudoinverse. The \textsf{Gappy-POD} approximation of the nonlinear vector is given as
	\begin{equation}\label{eq:gappypod}
		\widetilde{\bx{f}}\left( \bx{x}, \bm{\mu} \right) = \bx{U} \left( \bx{P}^{\textsf{T}} \bx{U} \right)^{\dagger} \bx{P}^{\textsf{T}} \bx{f}\left( \bx{x}, \bm{\mu} \right),		
	\end{equation}
where the matrix $\bx{P}$ now has $m > \ell$ columns and we have $\bx{P}^{\textsf{T}} \bx{U} \in \R^{m \times \ell}$. In \cite{morPehDG18}, the authors discuss two different oversampling strategies. The first approach, called \textsf{Gappy-POD Eigenvector}, considers the optimization point of view -- newly added interpolation points are those that lead to the largest decrease of $\| \left( \bx{P}^{\textsf{T}} \bx{U} \right)^{\dagger} \|$. The second oversampling strategy, \textsf{Gappy-POD Clustering}, proposed in \cite{morPehDG18} can be viewed as \textsf{Gappy-POD} based on interpreting the QR decomposition as a clustering algorithm. The additional samples from $\ell + 1$ till $m$ are identified based on the mutual entropy of the columns. For a more elaborate discussion we refer to \cite{morPehDG18,Zhaetal}. In the next section, we aim to make use of \textsf{DEIM} and its variants, as well as \textsf{Gappy-POD} to select \emph{important} parameter samples from the parameter domain.

For the pseudocode and the details, we refer the reader to the work \cite{morDrmG16} for the \textsf{QDEIM} algorithm and the work \cite{morPehDG18} in case of the \textsf{KDEIM} and \textsf{Gappy-POD} algorithms.
\section{Proposed Subsampling Strategy for the Training Set}
\label{sec:trngset}
A representative training set $\Xi_{\text{train}}$ is crucial for RBM to obtain a ROM that satisfies the required tolerance. While a densely-sampled $\Xi_{\text{train}}$ is needed to accurately represent the parameter space, it incurs high computational cost. In contrast, a randomly sampled coarse training set may fail to capture all the variations of the solution over the parameter space and result in a ROM that fails to meet the tolerance. Therefore, a wisely sampled coarse training set is desired to make the greedy algorithm efficient while retaining the required accuracy of the ROM. We now discuss two observations that motivate our proposed approach for training set sampling.	
\subsection{Motivating Observations} \label{subsec:motivation}
We detail two observations that pertain to the greedy algorithm in the RBM, the \textsf{DEIM} algorithm as well as the QR pivoting. We shall see that these two observations have motivated us to develop a subsampling strategy for the RBM training set.
\subsubsection{Greedy Parameters, QR Pivots, and DEIM Interpolation Points} \label{subsubsec:motivation1}
Our first observation concerns the parameters $\bm{\mu}^{*}$ selected by the greedy algorithm. The second observation is their resemblance to the QR pivots and the \textsf{DEIM} interpolation points. 

From our experience, the greedy algorithm tends to repeatedly pick parameter samples from a small subset of the training set, especially for time-dependent problems. This same phenomenon has been reported in other existing works \cite{morGre12, morHaa17, morHaaO08, morMadS13, morNguRP09, morQuaRM11, morZhaetal15}. The solution (or output) vectors at these parameter values usually exhibit large variability. While  the greedy algorithm scans through all the parameter samples in the training set, a majority of those samples are never picked. The fact that a few parameters get repeatedly picked reveals that there are still unresolved modes and hence more POD modes, corresponding to the selected parameter, are needed to get a good approximation. These few parameters picked by the greedy algorithm, usually represent solutions that are less smooth and hence are more difficult to approximate. An example of this phenomenon occurs in fluid dynamics problems where the low viscosity solutions develop shock and need a large number of POD modes to approximate. We illustrate this observation through the standard greedy algorithm applied to the discretized 1-D viscous Burgers' equation with $n = 1000$. The detailed description of the model is presented in \Cref{sec:numerics}. The parameter considered is the viscosity. We use $100$ equispaced parameter samples from the domain $\mathcal{P} = [0.005 \, , \, 1]$ to form the training set $\Xi_{\text{train}}$. A ROM with error below the tolerance $\texttt{tol} = 10^{-6}$ is requested from the greedy algorithm.
In \cref{tab:ill1a}, we provide the parameters picked by the greedy algorithm at each iteration from the training set. Noticeably, among the $100$ parameter samples in the training set, only $6$ contribute to generating the basis $\bx{V}$ that approximates the solution manifold. Of these $6$ samples, the sample $\mu = 0.005$ is picked fourteen times. This is not surprising since this parameter corresponds to the solution vector with the smallest viscosity and is the most difficult to approximate. 

Next, we make an important connection between the parameters selected by the greedy algorithm and the pivots obtained through a pivoted QR decomposition of the transpose of the output snapshot matrix defined in \cref{eq:y_sshot}.

For the same viscous Burgers' equation, we collect the snapshots of the scalar-valued outputs $\bx{y}$ at all the parameters in $\Xi_{\textnormal{train}}$ into a snapshot matrix given by
	\begin{equation}\label{eq:y_sshot}
		\mathbf{\mathfrak{Y}} := \begin{bmatrix}
									[\bx{y}(\bm{x}(t^{0},\bm{\mu}_{1}), \bm{\mu}_{1})]^{\textsf{T}} & \cdots & [\bx{y}(\bm{x}(t^{K}, \bm{\mu}_{1}), \bm{\mu}_{1})]^{\textsf{T}} \\
									\vdots & \ddots & \vdots \\
									[\bx{y}(\bm{x}(t^{0},\bm{\mu}_{N_{\text{train}}}), \bm{\mu}_{N_{\text{train}}})]^{\textsf{T}} & \cdots & [\bx{y}(\bm{x}(t^{K}, \bm{\mu}_{N_{\text{train}}}), \bm{\mu}_{N_{\text{train}}})]^{\textsf{T}}
								\end{bmatrix},
	\end{equation} 
Each row of the matrix consists of the snapshots of the output at $K+1$ time instances corresponding to a given parameter. Consider first the well-known pivoted QR decomposition of a matrix $\bx{D}$ given by
	\begin{equation}\label{eq:pivqr}
		\bx{D} \Pi = \bx{Q} \bx{R} = 	\begin{bmatrix}
							\bx{R}_{11} & \bx{R}_{12} \\
							0 & \bx{R}_{22}
						\end{bmatrix},
	\end{equation}
where $\bx{Q}$ is an orthogonal matrix and $\bx{R}$ is upper triangular. The pivots are given by the column permutation matrix $\Pi$.
\begin{table}[t!]
{\footnotesize
\captionsetup{position=top}
\caption{Greedy parameters and QR Pivots for the Burgers' equation.}
\label{tab:ill1}
\begin{center}
\subfloat[Greedy parameters picked by RBM.]{\label{tab:ill1a}
\begin{tabular}{|c|c|c|c|c|c|c|} \hline
Parameter & 0.005 & 0.0151 & 0.0251 & 0.0352 & 0.0553 & 1 \\ \hline
Repetitions & 14 & 1 & 1 & 1 & 1 & 1\\ \hline
\end{tabular}}
\\
\subfloat[First $10$ pivots for the QR decomposition of the transposed true output snapshot matrix $\mathbf{\mathfrak{Y}}$.]{\label{tab:ill1b}
\begin{tabular}{|c|c|c|c|c|c|c|c|c|c|c|} \hline
%Pivots & 0.005 & 0.0151 & 0.0251 & 0.0352 & 0.0653 & 1 \\ \hline
Pivots & 0.005 & 0.0151 & 0.0251 & 0.0352 & 0.0553 & 0.1055 & 0.1859 & 0.3166 & 0.7487 & 1 \\ \hline
\end{tabular}}
\end{center}
}
\end{table}
\setlength\fheight{1.5cm}
\setlength\fwidth{0.99\columnwidth}		
\begin{figure}[t!]
\centering
\begin{tikzpicture}

\begin{axis}[%
width=\fwidth,
height=\fheight,
at={(0\fwidth,0\fheight)},
scale only axis,
xmin=0.005,
xmax=1,
xticklabel style={
	/pgf/number format/fixed,
	/pgf/number format/precision=3
},
scaled x ticks=true,
xtick = {0.005, 0.1,0.2,0.3,0.4,0.5,0.6,0.7,0.8,0.9,1},
ymin=0,
ymax=2.22044604925031e-16,
ymajorticks=false,
axis y line=none,
axis x line*=bottom,
axis background/.style={fill=white},
%axis x line*=bottom,
%axis y line*=left,
legend columns=2,
legend style={only marks, at={(0.5,0.8)},anchor=north, font = \small,column sep=0.5cm}
]
\addplot [only marks, draw=red!80, thick, mark size=1.6pt, mark=o, mark options={}]
  table[row sep=crcr]{%
0.005	0\\
0.748737373737374	0\\
0.105505050505051	0\\
0.0251010101010101	0\\
0.316565656565657	0\\
1	0\\
0.0150505050505051	0\\
0.0552525252525252	0\\
0.185909090909091	0\\
0.0351515151515151	0\\
};
\addlegendentry{QR Pivots}

\addplot [color=black, draw=none, line width = 0.7pt, mark size=1.5pt, mark=x, mark options={solid, black, fill = black}]
  table[row sep=crcr]{%
0.005	0\\
0.0150505050505051	0\\
0.0251010101010101	0\\
0.0351515151515151	0\\
0.0552525252525252	0\\
1	0\\
};
\addlegendentry{Greedy Points}

\end{axis}
\end{tikzpicture}%
\caption{Greedy parameters for the Burgers' equation and QR pivots of the true output snapshots matrix $\mathbf{\mathfrak{Y}}$.}
\label{fig:ill1}
\end{figure}
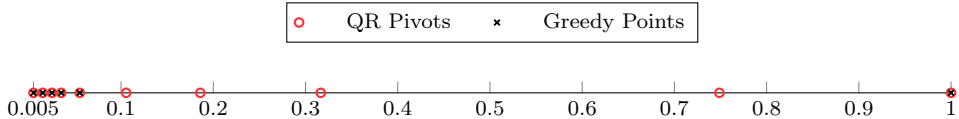
We apply the QR decomposition to $\mathbf{\mathfrak{Y}}^{\textsf{T}}$ and identify the pivots. A comparison of the parameters corresponding to the first ten pivots and the parameters selected in the greedy algorithm is shown in \cref{tab:ill1} and \cref{fig:ill1}. Of the ten pivots, six are identical with the greedy parameters. This close connection between the pivots of the QR decomposition and the greedy parameters chosen in the RBM has been, to the best of our knowledge, discussed only in \cite{morAntCF18, morNar20}. It is based on the interpretation of the QR decomposition as a greedy column selection procedure. Note that the application of a QR decomposition assumes the existence of the FOM solution for all the parameters in the training set. In practice, we do not have this information. Instead, we propose to apply the pivoted QR decomposition to the transpose of an approximate output snapshot matrix, in order to identify \emph{important} parameters which can then be used to subsample the fine training set in the RBM.

In \cref{subsec:deim}, the usage of \textsf{DEIM} and \textsf{QDEIM} algorithms were discussed in the context of MOR. The \textsf{DEIM} algorithm uses a greedy sparse sampling of the left singular vector matrix ($\bx{U}$) of the snapshots to identify interpolation points, whereas \textsf{QDEIM} performs a QR decomposition with column pivoting on $\bx{U}^{\textsf{T}}$. This implicates a similar phenomenon as observed above: QR with pivoting could select points of importance on different demands.  It is also noticed that QR with pivoting connects the greedy algorithm with \textsf{DEIM} (\textsf{QDEIM}), which indicates that \textsf{DEIM} and \textsf{QDEIM} could also be used to select representative parameter samples if either is applied to the output snapshot matrix.

\cref{fig:ill1} shows that the pivots of the QR decomposition on $\mathbf{\mathfrak{Y}}^{\textsf{T}}$ gives similar points as those selected by the greedy algorithm. It then indicates that sample points selected by the greedy algorithm in a way are highly related to the interpolation points of \textsf{QDEIM}, if the same snapshot matrix is considered by both the greedy algorithm and \textsf{QDEIM}, which is, in our case, the output snapshot matrix $\mathbf{\mathfrak{Y}}$. By exploiting this interpretation, we propose to use \textsf{DEIM} or other variants of \textsf{DEIM} in order to adapt the training set during the greedy algorithm.

To further support and motivate our proposed scheme, in the next subsection, we show that \textsf{DEIM} also has the similar capability of identifying the most representative parameter samples for dynamics, as that exhibited by the greedy algorithm in the RBM.
\setlength\fheight{4cm}
\setlength\fwidth{4cm}		
\begin{figure}[t]
\centering
\subfloat[Singular value decay.]{\label{fig:deimilla}\input{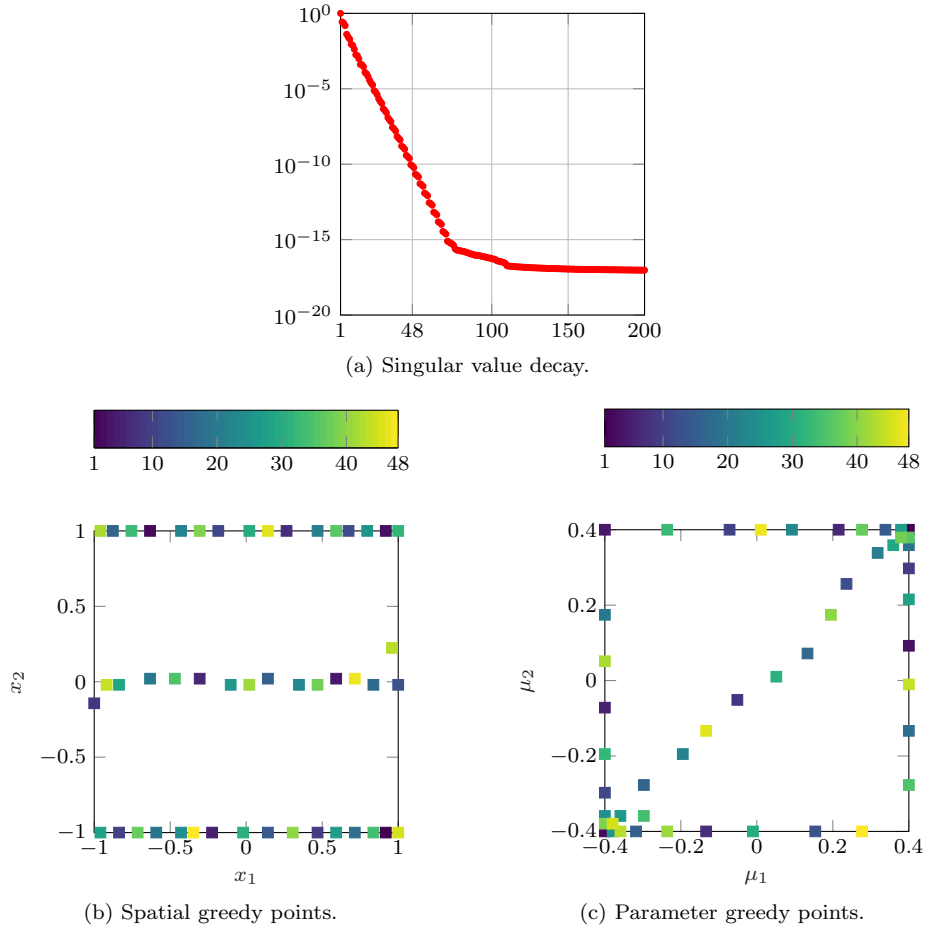}}\par	
\subfloat[Spatial greedy points.]{\label{fig:deimillb}% This file was created by matlab2tikz.
\begin{tikzpicture}

\begin{axis}[%
width=\fwidth,
height=\fheight,
at={(0\fwidth,0\fheight)},
scale only axis,
point meta min=1,
point meta max=48,
colormap/viridis,
colorbar horizontal,
colorbar sampled,
colorbar style={%
     xtick={1,10,20,30,40,48}, at={(0,1.4)},anchor=north west, samples = 48},
xmin=-1,
xmax=1,
xlabel style={font=\color{white!15!black}},
xlabel={$x_{1}$},
ymin=-1,
ymax=1,
ylabel style={font=\color{white!15!black}},
ylabel={$x_{2}$},
axis background/.style={fill=white},
colorbar
%  colormap={hot2}{
%  samples of colormap=(48 of hot2)},
%  colorbar right,
%  colorbar style={%
%     ytick={1, 10, 20, 30, 40},
%  } 
]
\addplot[scatter, only marks, mark=square*, mark size=2pt, scatter src=explicit, scatter/use mapped color={mark options={}, draw=mapped color, fill=mapped color}] table[row sep=crcr, meta=color]{%
x	y	color\\
0.918367346938775	-1	1\\
-0.63265306122449	1	2\\
0.918367346938775	1	3\\
0.591836734693878	0.0204081632653061	4\\
-0.224489795918367	-1	5\\
-0.306122448979592	0.0204081632653061	6\\
0.26530612244898	1	7\\
-1	-0.142857142857143	8\\
0.469387755102041	-1	9\\
-0.836734693877551	-1	10\\
-0.183673469387755	1	11\\
0.673469387755102	1	12\\
1	-0.0204081632653061	13\\
0.142857142857143	-1	14\\
0.142857142857143	0.0204081632653061	15\\
-0.877551020408163	1	16\\
0.836734693877551	-0.0204081632653061	17\\
0.714285714285714	-1	18\\
-0.591836734693878	-1	19\\
-0.63265306122449	0.0204081632653061	20\\
0.469387755102041	1	21\\
-0.428571428571429	1	22\\
0.346938775510204	-0.0204081632653061	23\\
-0.428571428571429	-1	24\\
0.795918367346939	1	25\\
-0.102040816326531	-0.0204081632653061	26\\
0.591836734693878	-1	27\\
-0.959183673469388	-1	28\\
-0.836734693877551	-0.0204081632653061	29\\
0.0204081632653061	1	30\\
-0.0204081632653061	-1	31\\
1	1	32\\
-0.755102040816326	1	33\\
0.836734693877551	-1	34\\
-0.469387755102041	0.0204081632653061	35\\
0.591836734693878	1	36\\
-0.714285714285714	-1	37\\
0.469387755102041	-0.0204081632653061	38\\
-0.306122448979592	1	39\\
0.306122448979592	-1	40\\
0.0204081632653061	-0.0204081632653061	41\\
-0.959183673469388	1	42\\
0.959183673469388	0.224489795918367	43\\
-0.918367346938776	-0.0204081632653061	44\\
1	-1	45\\
0.142857142857143	1	46\\
0.714285714285714	0.0204081632653061	47\\
-0.346938775510204	-1	48\\
};
\end{axis}
\end{tikzpicture}%
}\hfill
\subfloat[Parameter greedy points.]{\label{fig:deimillc}% This file was created by matlab2tikz.
\begin{tikzpicture}

\begin{axis}[%
width=\fwidth,
height=\fheight,
at={(0\fwidth,0\fheight)},
scale only axis,
point meta min=1,
point meta max=48,
colormap/viridis,
colorbar horizontal,
colorbar style={%
     xtick={1,10,20,30,40,48}, at={(0,1.4)},anchor=north west},
xmin=-0.4,
xmax=0.4,
xlabel style={font=\color{white!15!black}},
xlabel={$\mu_{1}$},
ymin=-0.4,
ymax=0.4,
ylabel style={font=\color{white!15!black}},
ylabel={$\mu_{2}$},
axis background/.style={fill=white},
colorbar
]
\addplot[scatter, only marks, mark=square*, mark size=2pt, scatter src=explicit, scatter/use mapped color={mark options={}, draw=mapped color, fill=mapped color}] table[row sep=crcr, meta=color]{%
x	y	color\\
0.4	0.4	1\\
-0.4	-0.4	2\\
0.4	0.0923076923076924	3\\
-0.4	0.4	4\\
-0.4	-0.0717948717948718	5\\
0.215384615384615	0.4	6\\
-0.133333333333333	-0.4	7\\
-0.0512820512820513	-0.0512820512820513	8\\
0.4	0.297435897435897	9\\
-0.0717948717948718	0.4	10\\
-0.4	-0.297435897435897	11\\
0.235897435897436	0.256410256410256	12\\
0.153846153846154	-0.4	13\\
0.338461538461538	0.4	14\\
-0.317948717948718	-0.4	15\\
-0.297435897435897	-0.276923076923077	16\\
0.133333333333333	0.0717948717948718	17\\
0.4	-0.133333333333333	18\\
0.4	0.358974358974359	19\\
-0.4	0.174358974358974	20\\
0.317948717948718	0.338461538461538	21\\
-0.194871794871795	-0.194871794871795	22\\
0.0923076923076924	0.4	23\\
-0.4	-0.358974358974359	24\\
-0.37948717948718	-0.4	25\\
-0.358974358974359	-0.358974358974359	26\\
0.379487179487179	0.4	27\\
0.4	0.215384615384615	28\\
0.358974358974359	0.358974358974359	29\\
-0.0102564102564103	-0.4	30\\
0.0512820512820513	0.0102564102564102	31\\
-0.4	-0.194871794871795	32\\
-0.235897435897436	0.4	33\\
-0.297435897435897	-0.358974358974359	34\\
0.4	-0.276923076923077	35\\
0.4	0.379487179487179	36\\
0.276923076923077	0.4	37\\
0.379487179487179	0.379487179487179	38\\
-0.4	-0.37948717948718	39\\
0.194871794871795	0.174358974358974	40\\
-0.235897435897436	-0.4	41\\
-0.4	0.0512820512820513	42\\
-0.358974358974359	-0.4	43\\
0.4	-0.0102564102564103	44\\
-0.37948717948718	-0.37948717948718	45\\
-0.133333333333333	-0.133333333333333	46\\
0.0102564102564102	0.4	47\\
0.276923076923077	-0.4	48\\
};
\end{axis}
\end{tikzpicture}%
}
\caption{Toy problem demonstrating anisotropic choice of interpolation points. The colourbars indicate the order of selection of the parameters. Points in the blue end of the spectrum were selected earlier while those in the green/yellow regions of the spectrum were picked later during the course of the algorithm.}
\label{fig:deimill}
\end{figure}
\subsubsection{DEIM and Parametric Anisotropy}
\label{subsubsec:motivation2}
For a function of two variables $\bx{f}(\bx{x}, \bm{\mu}) : \R^{n} \times \R^{p} \to \R^{n}$, the \textsf{DEIM} algorithm first identifies a linear subspace $\bx{U}$ and a small subset of points in the $\bx{x}$ variable, based on the snapshots matrix $\mathbf{\mathfrak{F}}$. One can analogously consider the mapping $\bx{f}(\bm{\mu}, \bx{x}) : \R^{p} \times \R^{n} \to \R^{n}$ through a transpose of $\bx{\mathfrak{F}}$. However, now the \textsf{DEIM} algorithm identifies a small subset of points in the $\bm{\mu}$ variable. We illustrate this on a toy example from \cite{Aan09}. 
Consider the following nonlinear, two parameter function:
	\begin{equation}
		\bx{f}(\bx{x}_{1}, \bx{x}_{2}, \bm{\mu}) = \dfrac{1 + \frac{\pi^{2}}{4}(\mu_{2} - \mu_{1} - (\mu_{1} + \mu_{2}) \bx{x}_{2})^{2} \sin^{2}(\frac{\pi}{2}(\bx{x}_{1} + 1) )}{1 + (\mu_{1} + \mu_{2}) \cos(\frac{\pi}{2}(\bx{x}_{1} + 1)) },
	\end{equation}
where $\bx{x} := [\bx{x}_{1}, \bx{x}_{2}] \in \R^{n \times 2}$ is the spatial variable obtained from the discretization of the two dimensional domain $\Omega := [-1 \,, 1] \times [-1 \,, 1]$ with 50 points in each spatial direction, resulting in $n = 2500$. The parameter $\bm{\mu} := (\mu_{1}, \mu_{2}) \in \R^{2}$ belongs to the domain $\mathcal{P} := [-0.4 \,, 0.4] \times [-0.4 \,, 0.4]$. 
We collect $1600$ snapshots of the function in the snapshots matrix $\mathbf{\mathfrak{F}} \in \R^{2500 \times 1600}$, based on uniform, equally spaced samples of the parameter. In \cref{fig:deimilla}, the singular values of the snapshot matrix $\mathbf{\mathfrak{F}}$ are plotted. The rapid decay clearly demonstrates the reducibility of this function. We apply a cut-off of $\epsilon_{\text{SVD}} = 10^{-10}$ for \cref{alg:deim} applied to both $\mathbf{\mathfrak{F}}$ and $\mathbf{\mathfrak{F}}^{\textsf{T}}$. The greedy interpolation points in \cref{fig:deimillb} are those corresponding to the indices in the row vector  $\bx{I}$, obtained from applying \cref{alg:deim} to $\mathbf{\mathfrak{F}}$. The greedy points in \cref{fig:deimillc} are those $\bm{\mu}$ corresponding to the indices stored in $\bx{I}$, obtained from applying \cref{alg:deim} to $\mathbf{\mathfrak{F}}^{\textsf{T}}$. The distribution of the points determined by \textsf{DEIM} for both the spatial and parameter variable have a characteristic structure. The number of interpolation points with SVD truncation tolerance $\epsilon_{\text{SVD}} = 10^{-10}$ was $48$, a mere $3\%$ of the total points. The spatial greedy points illustrate that while the variable $x_{1}$ is equally important over the entire range of $[-1 \,, 1]$, the $x_{2}$ variable has almost all its variation concentrated at $x_{2} \in \{-1, 0, 1\}$. For the parameter variable, the greedy algorithm picks most of the samples from the boundary of the domain and from the diagonal going from the lower left to the upper right. There is a dense concentration of points around the corners $(-0.4, -0.4)$ and $(0.4, 0.4)$. The choice of the greedy points is closely related to the structure of the function $\bx{f}$ being approximated. In most of the existing MOR literature, the \textsf{DEIM} algorithm has been used mainly as a tool to speed up evaluations of nonlinear or nonaffine (parametric) functions in a ROM. However, through the toy example, we have demonstrated its capability to expose the nature of parametric dependence of a function. As seen in \cref{fig:deimillc}, it is able to identify the regions in the parameter space where the function has large variations.

\subsection{Subsampling the Training Set}
Based on the observations in \Cref{subsubsec:motivation1,subsubsec:motivation2}, it is evident that a substantial computational effort can be saved in the offline stage of the RBM if we appropriately (optimally) sample the training set. The rationale for the proposed approach is the following: the standard greedy algorithm scans through the entire training set at each iteration and evaluates the error estimator at each parameter. This approach can incur significant computational cost for training sets with a large number of parameters. The proposed algorithm aims at picking out a small subset of the training set containing the most informative parameters. As will be demonstrated numerically, the parameters match closely to those chosen by the standard greedy algorithm.

Based on our observation of \cref{fig:ill1} in \Cref{subsubsec:motivation1} and \cref{fig:deimill} in \Cref{subsubsec:motivation2}, we propose to apply the pivoted QR decomposition and the \textsf{DEIM} algorithm (and its variants) to the snapshot matrix of the approximate output vector $\widetilde{\bx{y}}\left(\bx{z}(t^{k},\bm{\mu}), \bm{\mu}\right)$. More specifically, we consider the output snapshot matrix given by
	\begin{equation}\label{eq:yr_sshot}
		\widetilde{\mathbf{\mathfrak{Y}}} := \begin{bmatrix}
									[\widetilde{\bx{y}}(\bx{z}(t^{0},\bm{\mu}_{1}), \bm{\mu}_{1})]^{\textsf{T}} & \cdots & [\widetilde{\bx{y}}(\bx{z}(t^{K}, \bm{\mu}_{1}), \bm{\mu}_{1})]^{\textsf{T}} \\
									\vdots & \ddots & \vdots \\
									[\widetilde{\bx{y}}(\bx{z}(t^{0},\bm{\mu}_{N_{\text{train}}}), \bm{\mu}_{N_{\text{train}}})]^{\textsf{T}} & \cdots & [\widetilde{\bx{y}}(\bx{z}(t^{K}, \bm{\mu}_{N_{\text{train}}}), \bm{\mu}_{N_{\text{train}}})]^{\textsf{T}}
								\end{bmatrix} \in  \R^{N_{\text{train}} \times q N_{t}}
	\end{equation}
with each row containing the snapshots of the approximated output quantity at $K+1$ time instances corresponding to a given parameter sample. 
\begin{remark}
In case of stead-state systems with a single output we apply the proposed subsampling approach on the approximate state snapshots. For this, we define $\widetilde{\mathbf{\mathfrak{Y}}} := \widetilde{\bm{\mathfrak{X}}}^{\textsf{T}}$, $\widetilde{\bm{\mathfrak{X}}}$ being the snapshot matrix of the approximate state vector ($\widetilde{\bm{x}}(\bm{\mu}) = \bx{V} \bx{z}(\bm{\mu})$) such that
	\begin{displaymath}
		\widetilde{\bm{\mathfrak{X}}} := 	\left[\widetilde{\bm{x}}(\bm{\mu}_{1}), \cdots,  \widetilde{\bm{x}}(\bm{\mu}_{N_{\text{train}}})\right] \in \R^{n \times N_{\text{train}}}.
	\end{displaymath}
For steady-state systems with multiple outputs, we define $\widetilde{\mathbf{\mathfrak{Y}}} := \left[\widetilde{\bx{y}}(\bm{\mu}_{1}), \cdots, \widetilde{\bx{y}}(\bm{\mu}_{N_{\text{train}}}) \right]^{\textsf{T}} \in~{\R^{N_{\text{train}} \times q}}$.
\end{remark}
Note that $\widetilde{\mathbf{\mathfrak{Y}}}$ can be obtained from a coarse or low-fidelity ROM of the original system without doing FOM simulation at all the parameter samples. We propose two sampling strategies: (i) apply pivoted QR decomposition to $\widetilde{\mathbf{\mathfrak{Y}}}^{\textsf{T}}$, (ii) apply \textsf{DEIM} or its variants to $\widetilde{\mathbf{\mathfrak{Y}}}$, in order to identify the structure of the parametric dependence of the output variable.
\begin{algorithm}[t!]
\caption{Reduced Basis Method with Training Set Subsampling (Scheme~1)}
\label{alg:rbadapt1}
\begin{algorithmic}[1]
\REQUIRE{Training set ($\Xi_{\text{train}}^{f}$), tolerance for the ROM ($\texttt{tol}$), coarse tolerance ($\texttt{tol}^{c}$), maximum iteration ($n_{g}$).}
\ENSURE{$\bx{V}$.}
\STATE{Initialize $\bx{V} = [\,\,]$, $\epsilon = 1 + \texttt{tol}$, \texttt{iter} = 1.}
\STATE Select $\bm{\mu}^{*}_{(1)}$ randomly from $\Xi_{\text{train}}^{f}$.

\hrulefill \\
\textbf{Stage 1} 
\vspace{0.1cm}
\WHILE{$\epsilon > \texttt{tol}^{c}$}
		\STATE{Solve FOM \cref{eq:odess} at $\bm{\mu}^{*}_{(\texttt{iter})}$ and obtain snapshot matrix $\bm{\mathfrak{X}}$.}
		\STATE Set $\bar{\bm{\mathfrak{X}}} = \bm{\mathfrak{X}} - \bx{V} \bx{V}^{\textsf{T}} \bm{\mathfrak{X}}$.
		\STATE{Enrich $\bx{V}$ with $r_{\text{POD}}$ left singular vectors of  $\bar{\bm{\mathfrak{X}}}$ as
		
		$\bx{V} \leftarrow{} \texttt{orth}\left(\left[ \bx{V} \,\, \bx{U}_{\bar{\bm{\mathfrak{X}}}}(:\,,\,1 \, : \, r_{\text{POD}}) \right]\right)$.}
	\STATE{\texttt{iter} = \texttt{iter} + 1.}
	\STATE{Set $\bm{\mu}^{*}_{(\texttt{iter})} = \operatorname{arg} \operatorname{max} \limits_{\bm{\mu} \in \Xi_{\text{train}}^{f}} \Delta(\bm{\mu})$.}
	\STATE{$\epsilon = \Delta(\bm{\mu}^{*}_{(\texttt{iter})})$.}
\ENDWHILE
\\\hrulefill \\
\textbf{Stage 2}
\vspace{0.1cm}
\STATE Perform pivoted QR decomposition of $\widetilde{\mathbf{\mathfrak{Y}}}^{\textsf{T}}$, or apply \textsf{DEIM} or a \textsf{DEIM} variant to $\widetilde{\mathbf{\mathfrak{Y}}}$ and identify the indices $\mathbf{I}$ of the QR pivots or \textsf{DEIM} interpolation points.
\STATE Identify new training set $\Xi_{\textnormal{train}}$ using distribution of $\mathbf{I}$.
\WHILE{$\epsilon > \texttt{tol}\,\, \& \&\, \, \texttt{iter} \leq n_{g}$}
		\STATE{Solve FOM \cref{eq:odess} at $\bm{\mu}^{*}_{(\texttt{iter})}$ and obtain snapshot matrix $\bm{\mathfrak{X}}$.}
		\STATE Set $\bar{\bm{\mathfrak{X}}} = \bm{\mathfrak{X}} - \bx{V} \bx{V}^{\textsf{T}} \bm{\mathfrak{X}}$.
		\STATE{Enrich $\bx{V}$ with $r_{\text{POD}}$ left singular vectors of  $\bar{\bm{\mathfrak{X}}}$ as
		
		$\bx{V} \leftarrow{} \texttt{orth}\left(\left[ \bx{V} \,\, \bx{U}_{\bar{\bm{\mathfrak{X}}}}(:\,,\,1 \, : \, r_{\text{POD}}) \right]\right)$.}
	\STATE{\texttt{iter} = \texttt{iter} + 1.}
	\STATE{Set $\bm{\mu}^{*}_{(\texttt{iter})} = \operatorname{arg} \operatorname{max} \limits_{\bm{\mu} \in \Xi_{\text{train}}} \Delta(\bm{\mu})$.}
	
	\STATE{$\epsilon = \Delta(\bm{\mu}^{*}_{(\texttt{iter})})$.}
\ENDWHILE
\\\hrulefill
\end{algorithmic}
\end{algorithm}
Once the distribution of the interpolation points is identified, we can then adapt the training set for subsequent iterations of the greedy algorithm. We now outline the proposed approach and discuss different computational strategies.
\begin{algorithm}[t!]
\caption{Reduced Basis Method with Training Set Subsampling (Scheme~2)}
\label{alg:rbadapt2}
\begin{algorithmic}[1]
\REQUIRE{Training set ($\Xi_{\text{train}}^{f}$), tolerance for the ROM ($\texttt{tol}$), maximum iteration ($n_{g}$).}
\ENSURE{$\bx{V}$.}
\STATE{Initialize $\bx{V} = [\,\,]$, $\epsilon = 1 + \texttt{tol}$, \texttt{iter} = 1.}
\STATE Select $\bm{\mu}^{*}_{(1)}$ randomly from $\Xi_{\text{train}}^{f}$.

\hrulefill \\
\textbf{Stage 1} 
\vspace{0.1cm}
\WHILE{not terminated}
		\STATE{Solve FOM \cref{eq:odess} at $\bm{\mu}^{*}_{(\texttt{iter})}$ and obtain snapshot matrix $\bm{\mathfrak{X}}$.}
		\STATE Set $\bar{\bm{\mathfrak{X}}} = \bm{\mathfrak{X}} - \bx{V} \bx{V}^{\textsf{T}} \bm{\mathfrak{X}}$.
		\STATE{Enrich $\bx{V}$ with $r_{\text{POD}}$ left singular vectors of  $\bar{\bm{\mathfrak{X}}}$ as
		
		$\bx{V} \leftarrow{} \texttt{orth}\left(\left[ \bx{V} \,\, \bx{U}_{\bar{\bm{\mathfrak{X}}}}(:\,,\,1 \, : \, r_{\text{POD}}) \right]\right)$.}
	\STATE{\texttt{iter} = \texttt{iter} + 1.}
	\STATE{Set $\bm{\mu}^{*}_{(\texttt{iter})} = \operatorname{arg} \operatorname{max} \limits_{\bm{\mu} \in \Xi_{\text{train}}^{f}} \Delta(\bm{\mu})$.}
	
	\STATE{$\epsilon = \Delta(\bm{\mu}^{*}_{(\texttt{iter})})$.}
	\STATE Perform pivoted QR decomposition of $\widetilde{\mathbf{\mathfrak{Y}}}^{\textsf{T}}$, or apply \textsf{DEIM} or a \textsf{DEIM} variant to $\widetilde{\mathbf{\mathfrak{Y}}}$ and identify the indices $\mathbf{I}$ of the QR pivots or \textsf{DEIM} interpolation points.
	\IF{$\texttt{iter} \geq 2$}
		\STATE Check if $\operatorname{length}(\bx{I}_{\texttt{iter - 1}}) == \operatorname{length}(\bx{I}_{\texttt{iter}})$. \STATE If true $\mathbf{break}$ and proceed to \textbf{Stage 2}.
	\ENDIF
	
\ENDWHILE
\\\hrulefill \\
\textbf{Stage 2}
\vspace{0.1cm}
\STATE Identify new training set $\Xi_{\textnormal{train}}$ using distribution of $\mathbf{I}$.
\WHILE{$\epsilon > \texttt{tol}\,\, \& \&\, \, \texttt{iter} \leq n_{g}$}
		\STATE{Solve FOM \cref{eq:odess} at $\bm{\mu}^{*}_{(\texttt{iter})}$ and obtain snapshot matrix $\bm{\mathfrak{X}}$.}
		\STATE Set $\bar{\bm{\mathfrak{X}}} = \bm{\mathfrak{X}} - \bx{V} \bx{V}^{\textsf{T}} \bm{\mathfrak{X}}$.
		\STATE{Enrich $\bx{V}$ with $r_{\text{POD}}$ left singular vectors of  $\bar{\bm{\mathfrak{X}}}$ as
		
		$\bx{V} \leftarrow{} \texttt{orth}\left(\left[ \bx{V} \,\, \bx{U}_{\bar{\bm{\mathfrak{X}}}}(:\,,\,1 \, : \, r_{\text{POD}}) \right]\right)$.}
	\STATE{\texttt{iter} = \texttt{iter} + 1.}
	\STATE{Set $\bm{\mu}^{*}_{(\texttt{iter})} = \operatorname{arg} \operatorname{max} \limits_{\bm{\mu} \in \Xi_{\text{train}}} \Delta(\bm{\mu})$.}
	
	\STATE{$\epsilon = \Delta(\bm{\mu}^{*}_{(\texttt{iter})})$.}
\ENDWHILE
\\\hrulefill
\end{algorithmic}
\end{algorithm}

The proposed sampling procedure consists of two stages. The first stage is identical to the standard RBM procedure outlined in \cref{alg:rb}. A finely sampled training set $\Xi_{\text{train}}^{f}$ is used. We consider two different stopping criteria for the first stage --- (a) the first stage runs until the maximum estimated error is below a coarse tolerance denoted by $\texttt{tol}^{c} > \texttt{tol}$, where \texttt{tol} is the desired error tolerance for the final ROM, or (b) at two successive iterations, the number of \textsf{DEIM} interpolation points or QR pivots does not change. The value of $\texttt{tol}^{c}$ is user-defined and is of order $\mathcal{O}(1)$ in this work. Based on the two stopping criteria, two different schemes of training set subsampling are presented in \cref{alg:rbadapt1,alg:rbadapt2}, respectively. For both algorithms, we do not reset the value of \texttt{iter} at the end of Stage 1, so the final value of \texttt{iter} upon convergence for \cref{alg:rbadapt1,alg:rbadapt2} is the total number of iterations required by either algorithms to converge to the desired tolerance. It is worth pointing out that an \textit{a posteriori} error estimator \cite{morCheFB20} is used in both stages of the proposed algorithms so that the parameter sample picked at each iteration is the one at which an estimated {\textit{output}} error is the largest. Furthermore, the adaptive basis enrichment technique proposed in \cite{morCheFB20} is implemented and serves as a possible solution to the issue of repeated parameter selection. The number of POD modes corresponding to a selected parameter sample is adaptively decided: when the estimated error is large, a higher number of POD modes ($r_{\text{POD}}$) for the selected parameter are added; otherwise fewer POD modes are added. This reduces the chance of the same parameter sample being repeatedly chosen at subsequent iterations. This adaptive basis enrichment is implemented for both stages of our proposed method. We now discuss several practical computational strategies in connection with Steps 11 and 12 in \cref{alg:rbadapt1} and Steps 12, 13 and 16 in \cref{alg:rbadapt2}.

\begin{remark}
The active subspaces method (ASM) is an approach for parameter space reduction that has been recently applied in the context of MOR \cite{morConDW14,morTezBR18,morRozHS20} mainly for the case of scalar valued outputs. The ASM identifies a set of important directions in the parameter space onto which the parameter vectors are projected. It does this by means of Monte Carlo sampling of the gradients (with respect to the parameter) of the scalar-valued output quantity at a selection of parameter samples. The \emph{active subspaces} are the eigenspaces of the (truncated) covariance matrix of the gradients. Compared to ASM, our approach differs in two significant ways. Firstly, the proposed subsampling strategy is applicable to vector-valued output quantities. Secondly, ASM requires calculation of the gradient of the output. Moreover, the user still has to define the training set over which the gradient samples are acquired. Our proposed approach does not require the calculation of any additional quantity.
\end{remark}
\begin{remark}
Our proposed subsampling strategy occuring in Step 11 of \cref{alg:rbadapt1} and Step 10 of \cref{alg:rbadapt2} shares similarity with the column subset selection problem (CSSP) in the fields of numerical linear algebra and data mining \cite{Mah12}. For some general data matrix $\bx{D} \in \R^{N \times M}$, the CSSP aims to identify $h < M$ independent columns of the matrix $\bx{D}$ such that the residual $\| \bx{D} - \bx{P}_{h} \bx{D} \|$ is minimized. Here, $\bx{P}_{h} = \bx{S} \bx{S}^{\dagger}$ is a projection matrix and $\bx{S} \in R^{N \times h}$ consists of the $h$ extracted columns from $\bx{D}$. A number of algorithms, both deterministic and randomized, have been proposed to solve the CSSP \cite{BouMD09,BroBP10,CivM12}. One popular approach is to apply some variant of the QR decomposition (column-pivoted, rank-revealing, hybrid, etc.) either to the data matrix \bx{D} or to the transpose of the (truncated) left (or right) singular matrix $\bx{U}$ (or $\bx{W}$) of $\bx{D}$. If we consider $\bx{D}$ as the approximate output snapshot matrix, then our proposed algorithm using pivoted QR (or \textsf{QDEIM}) can be seen as a special case of the CSSP.
\end{remark}

\subsubsection{Training Set Subsampling: Scheme 1}
	We describe the proposed approach for the first scheme detailed in \cref{alg:rbadapt1}. In the first stage, a low-fidelity ROM is built with a coarse tolerance $\texttt{tol}^{c}$, over a finely sampled training set $\Xi_{\textnormal{train}}^{f}$. The intuition is that a low-fidelity approximation is sufficient to discover the parametric dependence of the output variable. At the end of the first stage, \textsf{DEIM} (or one of its variants) is applied to $\widetilde{\mathbf{\mathfrak{Y}}}$ to identify the interpolation points or a pivoted QR decomposition of $\widetilde{\mathbf{\mathfrak{Y}}}^{\textsf{T}}$ is used to identify the pivots. Once the set of interpolation points (or pivots) $p_{i}$ is identified, we proceed to suitably subsample the finely sampled training set based on the distribution of the identified interpolation points or pivots. Different possibilities exist to achieve this. A simple approach is to consider the training set for the second stage $\Xi_{\textnormal{train}}$ populated only by the identified interpolation points or pivots. Consider the fine training set $\Xi_{\textnormal{train}}^{f} := \{\bm{\mu}_{1}, \bm{\mu}_{2}, \ldots, \bm{\mu}_{N_{\text{train}}^{f}}\}$ with the subscript denoting the index corresponding to the position of a parameter in the set. Let $\bx{I}$ be the vector whose elements are the indices of the chosen interpolation points or pivots. We define the subsampled training set $\Xi_{\textnormal{train}}$ as the one consisting of all those parameters $\bm{\mu}_{z}$ from $\Xi_{\textnormal{train}}^{f}$ such that their indices are present in $\bx{I}$, i.e., $\Xi_{\textnormal{train}} := \{\bm{\mu}_{z; z \in \bx{I}}\}$. If there are only a few interpolation points, this approach would lead to a rapid second stage of the algorithm. However, there may exist the pitfall that it may result in an \textit{overfit} by which we mean that the resulting ROM after Stage 2 satisfies the desired tolerance \texttt{tol} only over the subsampled training set and does not generalize to other parameters in the parameter domain. We illustrate this phenomenon in the numerical tests. Another possible approach is to define a training budget $m$ for the second stage and use an oversampling strategy like the \textsf{Gappy-POD} to ensure that the training set for the second stage consists of a total of $m$ parameter samples.
\subsubsection{Training Set Subsampling: Scheme 2}	
	The first scheme of our proposed training set adaptation method requires a user-defined coarse tolerance. Choosing such a tolerance is rather heuristic. For some problems, a rough approximation may be enough to suitably capture all the parametric dependences, whereas a finer approximation may be needed for others. Therefore, for the second scheme we define a heuristic criterion that leads automatically to the second stage. To achieve this, we apply \textsf{DEIM} approximation (or a variant of \textsf{DEIM}) to the matrix $\widetilde{\mathbf{\mathfrak{Y}}}$ or the pivoted QR factorization to $\widetilde{\mathbf{\mathfrak{Y}}}^{\textsf{T}}$,  at each iteration of the greedy algorithm. Whenever the \textsf{DEIM} interpolations or pivoted QR decompositions at two successive iterations turn out to be equal, we terminate the first stage. This can be easily calculated by comparing the number of interpolation indices or the pivot indices at two successive iterations. Then, the subsampling of the training set for the second stage is carried out with similar approaches as discussed above. Following this, the second stage is run, until the required tolerance \texttt{tol} is met. In \cref{alg:rbadapt2}, $\bx{I}_{\texttt{iter}}$ and $\bx{I}_{\texttt{iter - 1}}$ in Stage 1 refer to the vector containing the interpolation indices identified by \textsf{DEIM} or its variants, or the indices of the QR pivots at the current iteration and the previous iteration, respectively.
	
Each of the two schemes has its own benefits or shortcomings. For Scheme 1, the burden of choosing an appropriate coarse tolerance lies with the user. This is highly problem-dependent. In the limit that $\texttt{tol}^{c} = \texttt{tol}$, Scheme 1 is just the standard Greedy algorithm (\cref{alg:rb}) with fixed training set. If $\texttt{tol}^{c} \gg \texttt{tol}$, a very fast second stage can be ensured, leading to considerable speedup of the offline stage of the RBM. However, it is not known \textit{a priori} if the chosen coarse tolerance is good. Scheme 2, on the other hand, automatizes the switching from Stage 1 to Stage 2 of the subsampling strategy by considering a heuristic criterion for the subspace approximation of the snapshot matrix $\widetilde{\mathbf{\mathfrak{Y}}}$. But, this comes with the additional cost of performing the \textsf{DEIM} algorithm or the pivoted QR decomposition, repeatedly. The success of both schemes is also highly dependent on the strategy adopted to construct the subsampled training set. In the numerical tests, we shall consider two approaches. In the first approach, we consider as many parameters in the subsampled training set as the number of \textsf{DEIM} interpolation points or QR pivots. For the second approach, we fix the cardinality of the subsampled training set to be $m$ and then use oversampling strategies based on the \textsf{Gappy-POD} method to choose those $m$ parameters in a principled way.
\subsubsection{Complexity Analysis}
\label{subsec:compl_analysis}

The fine training set is used in Stage 1 of both \Cref{alg:rbadapt1,alg:rbadapt2}. However, the computational complexity is not high for \Cref{alg:rbadapt1}, since we use a coarse tolerance $\texttt{tol}^{c}$ in Stage 1, so that the greedy algorithm converges within much fewer steps than when using the user desired tolerance in Stage 2. The computational complexity will grow with the decrease of the coarse tolerance used in Stage 1. However, as we have observed, a moderate value of $\texttt{tol}^{c}$ is enough to figure out the parameter dependency of the output. The number of FOM simulations is indeed independent of the size of the fine training set, since the FOM simulation is implemented only at those ``selected" parameter samples.
The situation for \Cref{alg:rbadapt2} is different since it involves the DEIM or QR algorithms at each iteration (see Step 10) to compute the interpolation points. A fine training set will indeed increase the cost of this step. Nevertheless, one can readily use recent techniques based on randomized linear algebra (such as randomized SVD, randomized QR, etc.) to keep the costs under control, see \cite{morSai20,morDueG17}.

In this section, we roughly compare the computational complexity of the standard RBM without training set subsampling in \cref{alg:rb} to that of RBM using the two subsampling strategies introduced in \cref{alg:rbadapt1,alg:rbadapt2}. Our complexity analysis considers the worst-case costs involved with each algorithm and is meant as a simplified illustration of the benefits of subsampling the training set. For simplicity, we only count the dominant costs in each algorithm.

We begin by introducing some notation:
\begin{itemize}
	\item Let $\mathcal{C}_{\text{FOM}} := O(n^{2})\cdot N_{t}$ denote the cost of solving the FOM (for e.g., \cref{eq:odess}) using some iterative scheme (GMRES, etc.), at a fixed parameter to obtain the snapshots.
	\item Let $\mathcal{C}_{\text{SVD}} := O(n^{3})$ be the worst-case cost for performing the SVD.
	\item Let $\mathcal{C}_{\text{Err}} := O(s)$ be the cost associated with evaluating the error estimator $\Delta(\bm{\mu})$ for a fixed parameter.
	\item Note that the dimension of the snapshot matrix in \cref{eq:yr_sshot} is $N_{\text{train}}^{f} \times q N_{t}$. For large-scale systems, usually we have $N_{\text{train}}^{f} \leq n$, if using, e.g., a sparse grid sampling, and often $q N_{t} \leq n$. Therefore, the cost for the \textsf{DEIM} algorithm (or its variants) or the pivoted QR decomposition should be less than $O(n^{3})$. In general, we denote it as $\mathcal{C}_{\text{ss}} := O(n^{3})$. 
\end{itemize}
The cost associated with \cref{alg:rb}, viz., $\mathcal{C}_{I}$ is:
\begin{align*}
	\mathcal{C}_{\RomanNumeralCaps{1}} := \left( \mathcal{C}_{\text{FOM}} + \mathcal{C}_{\text{SVD}} + \mathcal{C}_{\text{Err}} \cdot N_{\textnormal{train}}^{f}  \right) \cdot N_{I}
\end{align*}
with $N_{\textnormal{train}}^{f}$ being the cardinality of the fine training set and $N_{I}$ the number of iterations taken by the greedy algorithm to converge to the desired tolerance. The cost incurred by \cref{alg:rbadapt1}, $\mathcal{C}_{\RomanNumeralCaps{2}}$, is divided between Stage 1 and Stage 2. This is given by:
\begin{align*}
\mathcal{C}_{\RomanNumeralCaps{2}} :=& \overbrace{\big( \mathcal{C}_{\text{FOM}} + \mathcal{C}_{\text{SVD}} + \mathcal{C}_{\text{Err}} \cdot N_{\textnormal{train}}^{f}  \big) \cdot N_{II,1}}^{\text{Stage 1 cost}} \\ &+ \underbrace{\mathcal{C}_{\text{ss}} + \big( \mathcal{C}_{\text{FOM}} + \mathcal{C}_{\text{SVD}} + \mathcal{C}_{\text{Err}} \cdot N_{\textnormal{train}} \big) \cdot N_{II,2}}_{\text{Stage 2 cost}}.
\end{align*}
In the above expression, $N_{II,1}, N_{II,2}$ are, respectively, the number of iterations of the greedy algorithm in Stage 1 and Stage 2 of \cref{alg:rbadapt1}. Typically, $N_{II,2} > N_{II,1}$ since we use a coarse tolerance for Stage 1. Moreover, $N_{\textnormal{train}}$ is the cardinality of the subsampled training set $\Xi_{\textnormal{train}}$. Finally, the computational cost of \cref{alg:rbadapt2}, $\mathcal{C}_{\RomanNumeralCaps{3}}$, is as follows:
\begin{align*}
\mathcal{C}_{\RomanNumeralCaps{3}} :=& \overbrace{\big( \mathcal{C}_{\text{FOM}} + \mathcal{C}_{\text{SVD}} + \mathcal{C}_{\text{Err}} \cdot N_{\textnormal{train}}^{f} + \mathcal{C}_{\text{ss}}  \big) \cdot N_{III,1}}^{\text{Stage 1 cost}} \\ &+ \underbrace{\big( \mathcal{C}_{\text{FOM}} + \mathcal{C}_{\text{SVD}} + \mathcal{C}_{\text{Err}} \cdot N_{\textnormal{train}} \big) \cdot N_{III,2}}_{\text{Stage 2 cost}}
\end{align*}
with $N_{III,1}, N_{III,2}$ being, respectively, the number of iterations of the greedy algorithm in Stage 1 and Stage 2 of \cref{alg:rbadapt2}. For \cref{alg:rbadapt1} or \cref{alg:rbadapt2} to be computationally better alternatives to \cref{alg:rb}, we should have $\mathcal{C}_{I} > \mathcal{C}_{II}$ and $\mathcal{C}_{I} > \mathcal{C}_{III}$. 
\paragraph{Cost benefit of \cref{alg:rbadapt1}:}

To compare the costs of \cref{alg:rb} and \cref{alg:rbadapt1}, we look at the expression $\mathcal{C}_{I} - \mathcal{C}_{II}$.
\begin{align*}
	\mathcal{C}_{I} - \mathcal{C}_{II} :=&\,  \mathcal{C}_{\text{FOM}} (N_{I} - N_{II,1} - N_{II,2})\\ 
	&+ \mathcal{C}_{\text{SVD}} (N_{I} - N_{II,1} - N_{II,2})\\
	&+ \mathcal{C}_{\text{Err}} (N_{\textnormal{train}}^{f} \cdot N_{I} - N_{\textnormal{train}}^{f} \cdot N_{II,1} - N_{\textnormal{train}} \cdot N_{II,2})\\
	&- \mathcal{C}_{\text{ss}}.
\end{align*}
By assuming $N_{I} \approx N_{II,1} + N_{II,2}$ we have for the above expression
\begin{align*}
\mathcal{C}_{I} - \mathcal{C}_{II} \approx&\, \mathcal{C}_{\text{Err}} (N_{\textnormal{train}}^{f} - N_{\textnormal{train}}) N_{II,2} - \mathcal{C}_{\text{ss}}.
\end{align*}
For $N_{\textnormal{train}}^{f} \gg N_{\textnormal{train}}$, the first term would dominate leading to $\mathcal{C}_{I} > \mathcal{C}_{II}$. Thus, we see that computational cost incurred by \cref{alg:rbadapt1} is less than that of \cref{alg:rb}. In our analysis, we have assumed that $N_{I} \approx N_{II,1} + N_{II,2}$. While this helps to simplify the analysis, it does not always hold. As we shall see in \cref{sec:numerics}, for the example of Burgers' equation this assumption is true, while it does not hold for the thermal block example. A more relaxed assumption is $N_{I} \gtrapprox N_{II,1} + N_{II,2}$ using which, it can still be seen that $C_{I} > C_{II}$.

\paragraph{Cost benefit of \cref{alg:rbadapt2}:}

We look at the difference in costs $\mathcal{C}_{I} - \mathcal{C}_{III}$.
\begin{align*}
\mathcal{C}_{I} - \mathcal{C}_{III} :=&\,  \mathcal{C}_{\text{FOM}} (N_{I} - N_{III,1} - N_{III,2})\\ 
&+ \mathcal{C}_{\text{SVD}} (N_{I} - N_{III,1} - N_{III,2})\\
&+ \mathcal{C}_{\text{Err}} (N_{\textnormal{train}}^{f} \cdot N_{I} - N_{\textnormal{train}}^{f} \cdot N_{III,1} - N_{\textnormal{train}} \cdot N_{III,2})\\
&- \mathcal{C}_{\text{ss}} N_{III,1}.
\end{align*}
By making the assumption $N_{I} \approx N_{III,1} + N_{III,2}$, we have for the above expression
\begin{align*}
\mathcal{C}_{I} - \mathcal{C}_{III} \approx&\, \mathcal{C}_{\text{Err}} (N_{\textnormal{train}}^{f} - N_{\textnormal{train}}) N_{III,2} - \mathcal{C}_{\text{ss}} N_{III,1}.
\end{align*}
Usually $N_{\textnormal{train}}^{f} \gg N_{\textnormal{train}}$ and moreover, $N_{III,1} < N_{III,2}$. Therefore, the first term would dominate, leading to $\mathcal{C}_{I} > \mathcal{C}_{III}$. Once again, we see that the RBM using the subsampling strategy (\cref{alg:rbadapt2}) incurs a smaller cost when compared to \cref{alg:rb}. In the case of \cref{alg:rbadapt2}, just like earlier, the assumption $N_{I} \approx N_{III,1} + N_{III,2}$ need not always hold. From the numerical tests in the next section, we actually have $N_{I} \gtrapprox N_{III,1} + N_{III,2}$. The numerical results in \cref{tab:burgers_v1,tab:burgers_v2} show that \cref{alg:rbadapt1,alg:rbadapt2} indeed achieve speedups that can only be secured when they possess less computational complexity than \cref{alg:rb}.

\begin{remark}[High-dimensional parameter spaces]
	Both \Cref{alg:rbadapt1,alg:rbadapt2}, can be used when the parameter space is high-dimensional. The cost in Stage 2 will not be affected much, since a small training set identified from Stage 1 is used. The main increase in cost is due to the need to solve additional ROMs and estimate the error in Stage 1 of the proposed algorithms (Step 8 in \Cref{alg:rbadapt1,alg:rbadapt2}). If sparse grid sampling and a cheap error estimator are used, the cost will not increase fast. Actually, we can go a step further and make use of cheaply computable surrogate models of the error estimator as done in \cite{morCheFB19}. In that work, we considered a radial basis surrogate for the error estimator that is adaptively updated during the greedy algorithm. We only evaluate the actual error estimator over a few parameter samples. We then use this data to form a surrogate model, which can be used to evaluate the error for different parameter samples in the fine training set in Stage 1. 
\end{remark}
\begin{remark}[Role of output quantity of interest]
	In many cases, there is often the requirement, based on the application, to have a good approximation for the entire state vector. However, in this work we have specifically focussed on a \emph{goal-oriented} approach. For several applications, like in control systems or fluid dynamics, only a small number of state variables may be of interest. By focussing on those states alone, the resulting ROM dimension can be considerably lowered when compared to the case where the entire state needs to be well approximated. Also, it is indeed true that different output QoIs may have different influences on the parameter samples chosen. However, if some QoIs give rise to quantities with similar emphasis, then they may result in similar parameter samples being chosen. For example, a QoI defined by the mean of the solution over the spatial domain and a QoI defined by the sum of the solution over the spatial domain should result in similar samples. But a QoI defined by the mean of the solution over space probably gives different samples from the QoI defined by the maximum of the state over the spatial domain.
\end{remark}
\section{Numerical Results}
\label{sec:numerics}
We test the proposed adaptive training set subsampling algorithm on two examples. They are:
	\begin{enumerate}
	\item Viscous Burgers' equation with one parameter,
	\item thermal block with four parameters.
	\end{enumerate}
The first example is a nonlinear system, while the other is linear. All numerical tests are carried out in $\operatorname{MATLAB}\textsuperscript{\textregistered}$ 8.5.0.197613 (R2015a) on a laptop with $\operatorname{Intel} \textsuperscript{\textregistered} \operatorname{Core} \textsuperscript{\texttrademark}$ i5-7200U @ 2.5 GHz and 8 GB of RAM. Next, we describe the metrics used in all the numerical tests:
\begin{itemize}
\item The results of the proposed algorithms (\cref{alg:rbadapt1,,alg:rbadapt2}) are compared against a standard implementation of the POD-Greedy algorithm (\cref{alg:rb}) with a fixed training set. The implementation adopts Galerkin projection.
The number of POD modes $r_{\text{POD}},\, r_{\text{EI}} $ to enrich the RB and \textsf{DEIM} bases is determined at each iteration based on the adaptive approach proposed in \cite{morCheFB20}. We have also used the primal-dual error estimator proposed by the authors of \cite{morCheFB20} for our implementation of \cref{alg:rb,alg:rbadapt1,alg:rbadapt2}. The dual RB basis required for the error estimator is generated separately.
\item The fixed training set used for \cref{alg:rb} and the initial fine training set used for \cref{alg:rbadapt1,,alg:rbadapt2} are the same.
\item For \cref{alg:rbadapt1,,alg:rbadapt2}, we apply (i) the pivoted QR decomposition to the transposed approximate output snapshot matrix $\widetilde{\mathbf{\mathfrak{Y}}}^{\textsf{T}}$ and, (ii) the \textsf{DEIM} variants on the approximate output snapshot matrix $\widetilde{\mathbf{\mathfrak{Y}}}$ as two approaches to subsample the fine training set.
\item The cut-off criterion to determine the number of pivots ($h$) for the pivoted QR decomposition in \cref{eq:pivqr} of the approximate output snapshot matrix in \cref{eq:yr_sshot} is based on the magnitude of the diagonal elements in the upper triangular matrix $\bx{R}$, i.e., we set $h = q$ based on the the smallest $q$ such that $|\bx{R}(q+1,q+1)|/|\bx{R}(1,1)| < \epsilon_{\text{QR}}$, with $q \in \{1, 2, \ldots, \min(N_{\text{train}}, N_{t})\}$. The pivoted QR decomposition can effectively identify the rank of a matrix with a small diagonal $\bx{R}(q+1,q+1)$. Although there are cases when the column pivoted QR decomposition fails, they are rare in practice \cite{Cha87}.  A more robust rank-revealing QR factorization \cite{Cha87} can also be straightforwardly applied to our proposed subsampling algorithm.  However, in this work, we simply use column pivoted QR. The intrinsic $\operatorname{MATLAB}\textsuperscript{\textregistered}$ command \texttt{qr} is used with the options $\texttt{vector}$ enabled, i.e., we call $[\bx{Q}, \bx{R}, p_{\text{qr}}] = \texttt{qr}(\widetilde{\mathbf{\mathfrak{Y}}}^{\textsf{T}},\,\,\textnormal{\textquotesingle}\texttt{vector}\textnormal{\textquotesingle})$. It returns the pivot indices $p_{\text{qr}}$ as a vector, from which we select the first $h$ as our subsampling indices, i.e., $\bx{I} = p_{\text{qr}}(1:h)$. Here, $\bx{I}$ is the vector whose elements are the indices of the QR pivots and it lets us choose the subsampled training set $\Xi_{\text{train}}$ for Stage 2 of our proposed method, based on the fine training set $\Xi_{\textnormal{train}}^{f}$ from Stage 1.
\item Our implementation of the \texttt{k-means} algorithm for \textsf{KDEIM} is based on the intrinsic $\operatorname{MATLAB}\textsuperscript{\textregistered}$ function \texttt{kmeans}. We use five different initializations and pick the best configuration among the five.
\item The maximum true error over the test set is defined as 
	\begin{displaymath}
		\epsilon_{\text{t}}^{\text{max}} := \max \limits_{\bm{\mu} \in \Xi_{\textnormal{test}}}  \left( \frac{1}{K+1} \sum_{k = 0}^{K} \| \bx{y}\left(\bx{x}(t^{k},\bm{\mu})\right) - \widetilde{\bx{y}}\left(\bx{z}(t^{k},\bm{\mu})\right) \| \right).
	\end{displaymath}
\item Reported runtimes for all the algorithms are obtained by considering the median value of five independent runs.
\item The quantity Iterations reported in \cref{tab:burgers_v1,tab:burgers_v2,tab:thermal_v1,tab:thermal_v2,tab:thermal_oversample} refers to the total number of iterations (\texttt{iter}) of the corresponding greedy algorithm (\cref{alg:rb,alg:rbadapt1,alg:rbadapt2}) to converge to the desired tolerance. 
\end{itemize}
\subsection{Viscous Burgers' Equation}
We consider the following viscous Burgers' equation defined on a 1-D domain $\Omega := [ 0 \,, 1]$:
	\begin{equation}
	\begin{aligned}
	 \frac{\partial x}{\partial t} + x \frac{\partial x}{\partial w} &= \mu \frac{\partial^{2} x}{\partial w^{2}} + s(x,t), \,\,\, x(0, t) = 0 \,\, \&\,\, \frac{\partial x(1,t)}{\partial w} = 0,\\
	 y &= x(1,t),
	 \end{aligned}
	\end{equation}
where $s(x,t)$ denotes a forcing term defined later. The domain is discretized using the finite difference method with step size $\Delta w = 0.001$. We make use of a second order central difference discretization for the diffusion term and a first-order upwind scheme for the convection term. The resulting FOM is of dimension $n = 1000$. For the time variable $t \in [0 \,, 2]$, we make use of a first order implicit-explicit (IMEX) scheme with the diffusion term discretized implicitly and the nonlinear convection term discretized explicitly. The step size is $\Delta t = 0.001$. The constant input to the system is set to be $s(x,t) \equiv 1$ and the initial condition is $\bx{x}_{0} := \mathbf{0} \in \R^{n}$. The parameter domain of the viscosity $\mu$ is $\mathcal{P} := [0.005, \, 1]$. The training set $\Xi_{\text{train}}$ consists of $100$ equally spaced samples in $\mathcal{P}$. For the RBM, we fix the tolerance to be $\texttt{tol} = 1\cdot10^{-6}$. To validate the ROM, we use a test set $\Xi_{\textnormal{test}}$ containing 300 randomly sampled parameters, different from those in the training set.
\subsubsection{Greedy Algorithm with Fixed Training Set}
We begin by applying the standard greedy algorithm (\cref{alg:rb}) to the discretized model of the Burgers' equation. The greedy algorithm requires $t_{\text{greedy}} = 505.47$ seconds and $19$ iterations to converge to the defined tolerance of $1\cdot10^{-6}$. The resulting ROM has RB dimension $r_{\text{POD}} = 32$ along with $r_{\text{EI}} = 33$ basis vectors for the \textsf{DEIM} projection matrix. The maximum true error over the test set is $\epsilon_{\text{t}}^{\text{max}} = 2.07\cdot10^{-8}$. Although the POD-Greedy algorithm with a fixed training set results in a ROM that meets the specified tolerance, its offline time is high and there is scope for improvement by considering a subsampled training set.
\subsubsection{Greedy Algorithm Schemes 1 and 2}
We apply \cref{alg:rbadapt1,,alg:rbadapt2} to the Burgers' equation, making use of both pivoted QR and the two \textsf{DEIM} variants (\textsf{QDEIM} and \textsf{KDEIM}) to identify the interpolation points $\bx{I}$. Further, we consider three different SVD and QR cut-off tolerances ($\epsilon_{\text{SVD}}$, $\epsilon_{\text{QR}}$) for the pivoted QR and \textsf{DEIM} variants $\{1\cdot10^{-4}, 1\cdot10^{-6}, 1\cdot10^{-8} \}$ to highlight the progressive refinement of the adapted training set in `difficult regions' of the parameter space. The results are summarized in \cref{tab:burgers_v1,,tab:burgers_v2} for \cref{alg:rbadapt1,,alg:rbadapt2}, respectively. The matrix $\widetilde{\mathbf{\mathfrak{Y}}} \in \R^{100 \times 81}$, for either algorithms, was assembled by collecting the snapshots of the output vector at every $25^{\text{th}}$ time step. The training set in Stage 2 of both algorithms consists of interpolation points identified by QR, \textsf{QDEIM} or \textsf{KDEIM}. As revealed in the results, this choice is sufficient to produce ROMs that meet the required tolerance over the test set. For this example, there is not a big difference between the results of the two algorithms. Both schemes produce ROMs of almost identical RB, \textsf{DEIM} basis sizes ($r_{\text{POD}},\,r_{\text{EI}}$) and result in nearly the same maximum error over the test set. 

We show the subsampled training sets resulting from \cref{alg:rbadapt2} using the pivoted QR, \textsf{QDEIM} and \textsf{KDEIM} variants, with different SVD, QR tolerances in \cref{fig:burgtrngset_minus4,fig:burgtrngset_minus6,fig:burgtrngset_minus8}. The black crosses denote those samples from the fine training set which were retained in Stage 2 of the algorithm. For the \textsf{QDEIM} variant, it is clear that the subsampled parameters are concentrated more around the lower viscosity regions of the parameter space. Thus, the method is able to successfully identify the \emph{physically more relevant} points. Moreover, the parameter samples identified by \textsf{QDEIM} are very close to the ones identified by the method using a pivoted QR decomposition. This is not surprising since the former determines the interpolation points through a pivoted QR decomposition of $\bx{U}^{\textsf{T}}$ ($\bx{U}$ is the left singular matrix of $\widetilde{\mathbf{\mathfrak{Y}}}$) whereas the latter applies the pivoted QR decomposition directly to $\widetilde{\mathbf{\mathfrak{Y}}}^{\textsf{T}}$. We also show the results of \textsf{KDEIM}, where the subsampled (selected) parameter samples and their corresponding clusters are presented. The subsampled points in this case are the centroids of the clusters. The clusters are smaller in size for the low viscosity regions, while they are comparatively larger in the high viscosity zone. The resulting subsampled training sets from \cref{alg:rbadapt1} display a similar trend.

For a given $\epsilon_{\text{SVD}}$ or $\epsilon_{\text{QR}}$, the runtimes for the \textsf{QDEIM} and \textsf{KDEIM} based training set adaptation are very close. Using the pivoted QR leads to a subsampled training set with one sample more than that generated by using \textsf{QDEIM} or \textsf{KDEIM}. This results in a marginally higher offline time for this method. One observation worth remarking is that, for some instances, using the \textsf{KDEIM} approach to identify the adapted training set leads to the greedy algorithm converging in fewer iterations. This is most likely due to the fact that the identified parameters in this case represent cluster centroids and are more representative of the average behaviour. This yields a more uniform approximation throughout the parameter domain, with each cluster average being well-represented. The \textsf{QDEIM} version, on the other hand, tends to identify points based on the SVD of the output snapshot matrix and tends to favour points away from the mean behaviour. We illustrate this in \cref{fig:err_plot_alg41} for the case of $\epsilon_{\text{SVD}} = 10^{-8}$ for \cref{alg:rbadapt1}. It is evident that while the subsampling strategy using \textsf{QDEIM} results in a smaller magnitude of the maximum error over the test set ($4.55\cdot10^{-8}$ vs. $6.88\cdot10^{-8}$), the approach using the \textsf{KDEIM}-based sampling leads to a more uniform distribution of the error over the test set.

\begin{landscape}
\centering
\begin{table}[t!]
\caption{Results of \cref{alg:rbadapt1} with varying $\epsilon_{\text{SVD}}$, $\epsilon_{\text{QR}}$ for Burgers' equation.}
\label{tab:burgers_v1}
\centering
\scriptsize
\begin{tabular}{|c|c|c|c|c|c|c|c|c|c|c|}
\hline
\multirow{3}{*}{Method}    & \multirow{3}{*}{Fixed} & \multicolumn{9}{c|}{Adapted}                                                                                  \\ \cline{3-11} 
                           &                        & \multicolumn{3}{c|}{$\epsilon_{\text{SVD}}, \epsilon_{\text{QR}} = 1\cdot10^{-4}$}          & \multicolumn{3}{c|}{$\epsilon_{\text{SVD}}, \epsilon_{\text{QR}} = 1\cdot10^{-6}$}           & \multicolumn{3}{c|}{$\epsilon_{\text{SVD}}, \epsilon_{\text{QR}} = 1\cdot10^{-8}$}          \\ \cline{3-11} 
                           &                        & \small{\textsf{QDEIM}}       			  & \small{\textsf{KDEIM}}			& \small{QR}
                           & \small{\textsf{QDEIM}} & \small{\textsf{KDEIM}}  & \small{QR}                & \small{\textsf{QDEIM}}
                           & \small{\textsf{KDEIM}} & \small{QR}          \\ \hline
$N_{\text{train}}$                    & 100                                 & 7                  & 7      & 8            & 11                  & 11	&		12                 & 15                 & 15					 & 15                \\ \hline
$\epsilon_{t}^{\text{max}}$ & $2.07 \cdot 10^{-8}$                & $4.22\cdot10^{-7}$ & $4.15\cdot10^{-7}$ & $3.43\cdot 10^{-8}$ &$1.87\cdot 10^{-8}$ & $1.87\cdot10^{-8}$ & $3.43\cdot 10^{-8}$  & $1.87\cdot10^{-8}$ & $1.84\cdot10^{-8}$ & $3.28\cdot 10^{-8}$\\ \hline
($r_{\text{POD}}, r_{\text{EI}}$)                 & (32,33)                             & (31,32)            & (30,31)   & (31,32)            & (31,32)             & (31,32)		& (31,32)            & (31,32)            & (31,32)	&(31,32)            \\ \hline
Iterations (\texttt{iter})                     & 19                                  & 17                 & 17     & 17            & 18                  & 18       & 17          & 18                 & 18           & 17      \\ \hline
Offline time (s)             & 505.47                              & 91.00              & 91.16       & 94.90       & 110.10              & 110.35    & 109.54         & 125.07             & 125.59        & 120.53             \\ \hline
Speedup                    & -                                   & 5.6                & 5.5           & 5.3     & 4.6                 & 4.6        & 4.6        & 4.0                & 4.0          & 4.2      \\ \hline
\end{tabular}
\end{table}
\hfill
\begin{table}[b!]
\caption{Results of \cref{alg:rbadapt2} with varying $\epsilon_{\text{SVD}}$, $\epsilon_{\text{QR}}$ for Burgers' equation.}
\label{tab:burgers_v2}
\centering
\scriptsize
\begin{tabular}{|c|c|c|c|c|c|c|c|c|c|c|}
\hline
\multirow{3}{*}{Method}    & \multirow{3}{*}{Fixed} & \multicolumn{9}{c|}{Adapted}                                                                                    \\ \cline{3-11} 
                           &                                     & \multicolumn{3}{c|}{$\epsilon_{\text{SVD}}, \epsilon_{\text{QR}} = 1\cdot10^{-4}$}          & \multicolumn{3}{c|}{$\epsilon_{\text{SVD}}, \epsilon_{\text{QR}} = 1\cdot10^{-6}$}           & \multicolumn{3}{c|}{$\epsilon_{\text{SVD}}, \epsilon_{\text{QR}} = 1\cdot10^{-8}$}          \\ \cline{3-11} 
                           &                                     & \small{\textsf{QDEIM}}              & \small{\textsf{KDEIM}}     & \small{QR}         & \small{\textsf{QDEIM}}               & \small{\textsf{KDEIM}}     & \small{QR}         & \small{\textsf{QDEIM}}              & \small{\textsf{KDEIM}}       & \small{QR}       \\ \hline
$N_{\text{train}}$                    & 100                                 & 7                  & 7       & 8          & 11                  & 11	 & 12                 & 15                 & 15             & 15    \\ \hline
$\epsilon_{t}^{\text{max}}$ & $2.07 \cdot 10^{-8}$                & $4.22\cdot10^{-7}$ & $3.43\cdot10^{-7}$ & $3.43\cdot10^{-8}$ & $1.87\cdot 10^{-8}$ & $5.96\cdot10^{-8}$ &  $3.43\cdot10^{-8}$ & $1.87\cdot10^{-8}$ & $2.59\cdot10^{-8}$  & $3.28\cdot10^{-8}$\\ \hline
($r_{\text{POD}}, r_{\text{EI}}$)                 & (32,33)                             & (31,32)            & (30,31)     & (31,32)       & (31,32)             & (31,32)    & (31,32)        & (31,32)            & (30,31)       & (31,32)      \\ \hline
Iterations (\texttt{iter})                      & 19                                  & 17                 & 16         & 17        & 18                  & 18       & 17          & 18                 & 16               & 17  \\ \hline
Offline time (s)              & 505.47                              & 90.78              & 88.52         & 94.58     & 110.12              & 110.10      & 109.55       & 145.22             & 135.82             & 120.27 \\ \hline
Speedup                    & -                                   & 5.6                & 5.7              & 5.3   & 4.6                 & 4.6         & 4.6       & 3.5                & 3.7        & 4.2        \\ \hline
\end{tabular}
\end{table}
\end{landscape}
On average, for a given SVD tolerance,  \cref{alg:rbadapt1} is faster than \cref{alg:rbadapt2}. This can be attributed to the fact that for the latter, the \textsf{DEIM} variant or the pivoted QR decomposition needs to be performed repeatedly to check the criterion in Step 12 of \cref{alg:rbadapt2}. Since this involves performing an SVD, the associated costs are higher. The proposed subsampling algorithms result in a noticeable speedup of the POD-Greedy algorithm. For \cref{alg:rbadapt1}, the maximum achieved speedup was $5.6$ for $\epsilon_{\text{SVD}} = 1\cdot10^{-4}$ using \textsf{QDEIM}. However, the least speedup noticed was $4.0$ for $\epsilon_{\text{SVD}} = 1\cdot10^{-8}$ using \textsf{QDEIM}. Also, for \cref{alg:rbadapt2} the maximum achieved speedup was $5.7$ for $\epsilon_{\text{SVD}} = 1\cdot10^{-4}$ using \textsf{KDEIM} while the minimum speedup was $3.5$ for $\epsilon_{\text{SVD}} = 1\cdot10^{-8}$ using \textsf{QDEIM}. Finally, we see from \cref{tab:burgers_v1,tab:burgers_v2} that the number of iterations to converge for \cref{alg:rb}, \cref{alg:rbadapt1,alg:rbadapt2} are nearly the same. Thus, our assumption in the analysis from \cref{subsec:compl_analysis} that $N_{I} \approx N_{II,1} + N_{II,2}$ and $N_{I} \approx N_{III,1} + N_{III,2}$ holds true. 
\setlength\fheight{4cm}
\setlength\fwidth{0.9\columnwidth}
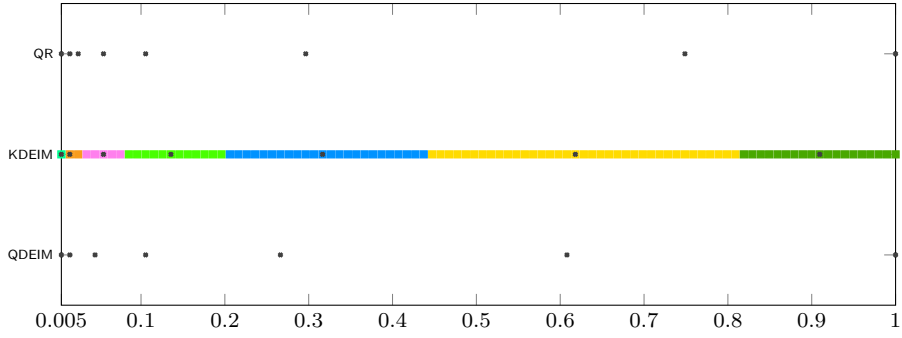
\begin{figure}[tp!]
	\centering
%\definecolor{mycolor1}{rgb}{1.00000,0.85714,0.00000}%
%\definecolor{mycolor2}{rgb}{0.28571,1.00000,0.00000}%
%\definecolor{mycolor3}{rgb}{0.00000,1.00000,0.57143}%
%\definecolor{mycolor4}{rgb}{0.00000,0.57143,1.00000}%
%\definecolor{mycolor5}{rgb}{0.28571,0.00000,1.00000}%
%\definecolor{mycolor6}{rgb}{1.00000,0.00000,0.85714}%
\definecolor{mycolor1}{rgb}{1.00000,0.85714,0.00000}%
\definecolor{mycolor2}{rgb}{0.28571,1.00000,0.00000}%
\definecolor{mycolor3}{rgb}{0.00000,1.00000,0.57143}%
\definecolor{mycolor4}{rgb}{0.00000,0.57143,1.00000}%
%\definecolor{mycolor5}{rgb}{0.28571,0.00000,1.00000}%
\definecolor{mycolor5}{rgb}{0.28571,0.65714,0.0000}%
\definecolor{mycolor6}{rgb}{1.00000,0.00000,0.85714}%
\begin{tikzpicture}

\begin{axis}[%
width=\fwidth,
height=\fheight,
at={(0\fwidth,0\fheight)},
scale only axis,
xmin=0.005,
xmax=1,
xticklabel style={
	/pgf/number format/fixed,
	/pgf/number format/precision=3
},
yticklabel style={font=\tiny},
xtick={0.005,   0.1,   0.2,   0.3,   0.4,   0.5,   0.6,   0.7,   0.8,   0.9,     1},
ymin=0.5,
ymax=3.5,
ytick={1,2,3},
yticklabels={{\textsf{QDEIM}},{\textsf{KDEIM}},{\textsf{QR}}},
axis background/.style={fill=white}
]
% QDEIM
\addplot [color=darkgray, line width=1.0pt, draw=none,mark size=1.0pt, mark=x, mark options={solid, darkgray}]
  table[row sep=crcr]{%
0.005	1\\
0.0150505050505051	1\\
0.0452020202020202	1\\
0.105505050505051	1\\
1	1\\
0.266313131313131	1\\
0.608030303030303	1\\
};
% KDEIM
\addplot [color=mycolor1, line width=0.5pt, draw=none, mark size=1.3pt, mark=square*, mark options={solid, mycolor1}]
  table[row sep=crcr]{%
0.447222222222222	2\\
0.457272727272727	2\\
0.467323232323232	2\\
0.477373737373737	2\\
0.487424242424242	2\\
0.497474747474748	2\\
0.507525252525252	2\\
0.517575757575758	2\\
0.527626262626263	2\\
0.537676767676768	2\\
0.547727272727273	2\\
0.557777777777778	2\\
0.567828282828283	2\\
0.577878787878788	2\\
0.587929292929293	2\\
0.597979797979798	2\\
0.608030303030303	2\\
0.618080808080808	2\\
0.628131313131313	2\\
0.638181818181818	2\\
0.648232323232323	2\\
0.658282828282828	2\\
0.668333333333333	2\\
0.678383838383838	2\\
0.688434343434343	2\\
0.698484848484848	2\\
0.708535353535354	2\\
0.718585858585859	2\\
0.728636363636364	2\\
0.738686868686869	2\\
0.748737373737374	2\\
0.758787878787879	2\\
0.768838383838384	2\\
0.778888888888889	2\\
0.788939393939394	2\\
0.798989898989899	2\\
0.809040404040404	2\\
};
\addplot [color=yellow!55!red, line width=0.5pt, draw=none, mark size=1.3pt, mark=square*, mark options={solid, yellow!55!red}]
  table[row sep=crcr]{%
0.0150505050505051	2\\
0.0251010101010101	2\\
};
\addplot [color=mycolor3, line width=0.5pt, draw=none, mark size=1.3pt, mark=square*, mark options={solid, mycolor3}]
  table[row sep=crcr]{%
0.005	2\\
};
\addplot [color=mycolor2, line width=0.5pt, draw=none, mark size=1.3pt, mark=square*, mark options={solid,mycolor2}]
  table[row sep=crcr]{%
0.0854040404040404	2\\
0.0954545454545455	2\\
0.105505050505051	2\\
0.115555555555556	2\\
0.125606060606061	2\\
0.135656565656566	2\\
0.145707070707071	2\\
0.155757575757576	2\\
0.165808080808081	2\\
0.175858585858586	2\\
0.185909090909091	2\\
0.195959595959596	2\\
};
\addplot [color=mycolor5, line width=0.5pt, draw=none, mark size=1.3pt, mark=square*, mark options={solid, mycolor5}]
  table[row sep=crcr]{%
0.819090909090909	2\\
0.829141414141414	2\\
0.839191919191919	2\\
0.849242424242424	2\\
0.859292929292929	2\\
0.869343434343434	2\\
0.879393939393939	2\\
0.889444444444444	2\\
0.899494949494949	2\\
0.909545454545455	2\\
0.91959595959596	2\\
0.929646464646465	2\\
0.93969696969697	2\\
0.949747474747475	2\\
0.95979797979798	2\\
0.969848484848485	2\\
0.97989898989899	2\\
0.989949494949495	2\\
1	2\\
};
\addplot [color=mycolor6!50, line width=0.5pt, draw=none, mark size=1.3pt, mark=square*, mark options={solid, mycolor6!50}]
  table[row sep=crcr]{%
0.0351515151515151	2\\
0.0452020202020202	2\\
0.0552525252525252	2\\
0.0653030303030303	2\\
0.0753535353535354	2\\
};
\addplot [color=mycolor4, line width=0.5pt, draw=none, mark size=1.3pt, mark=square*, mark options={solid, mycolor4}]
  table[row sep=crcr]{%
0.206010101010101	2\\
0.216060606060606	2\\
0.226111111111111	2\\
0.236161616161616	2\\
0.246212121212121	2\\
0.256262626262626	2\\
0.266313131313131	2\\
0.276363636363636	2\\
0.286414141414141	2\\
0.296464646464646	2\\
0.306515151515152	2\\
0.316565656565657	2\\
0.326616161616162	2\\
0.336666666666667	2\\
0.346717171717172	2\\
0.356767676767677	2\\
0.366818181818182	2\\
0.376868686868687	2\\
0.386919191919192	2\\
0.396969696969697	2\\
0.407020202020202	2\\
0.417070707070707	2\\
0.427121212121212	2\\
0.437171717171717	2\\
};
\addplot [color=darkgray, line width=1.0pt, draw=none, mark size=1.0pt, mark=x, mark options={solid, darkgray}]
  table[row sep=crcr]{%
0.618080808080808	2\\
0.0150505050505051	2\\
0.005	2\\
0.135656565656566	2\\
0.909545454545455	2\\
0.0552525252525252	2\\
0.316565656565657	2\\
};
% QR
\addplot [color=darkgray, line width=1.0pt,mark size=1.0pt, draw=none, mark=x, mark options={solid, darkgray}]
  table[row sep=crcr]{%
0.005	3\\
0.748737373737374	3\\
0.105505050505051	3\\
0.0251010101010101	3\\
0.296464646464646	3\\
1	3\\
0.0150505050505051	3\\
0.0552525252525252	3\\
};
\end{axis}
\end{tikzpicture}%
	\caption{\cref{alg:rbadapt2} for the Burgers' Equation with SVD, QR tolerance $\epsilon_{\text{SVD}}, \epsilon_{\text{QR}} = 1\cdot10^{-4}$. The crossmarks denote the parameters in the subsampled training set. For \textsf{KDEIM} each colour represents one cluster; the centroids of each of the clusters make up the subsampled training set.}
	\label{fig:burgtrngset_minus4}
\end{figure}
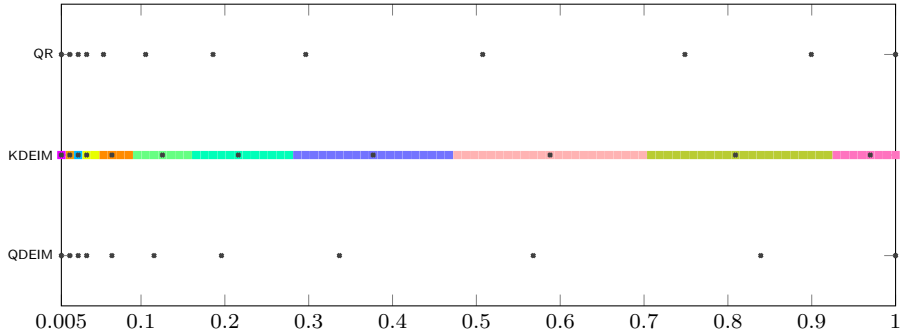
\begin{figure}[tp!]
	\centering
\definecolor{mycolor1}{rgb}{1.00000,0.54545,0.00000}%
\definecolor{mycolor2}{rgb}{0.90909,1.00000,0.00000}%
\definecolor{mycolor3}{rgb}{0.00000,1.00000,0.18182}%
\definecolor{mycolor4}{rgb}{0.00000,1.00000,0.72727}%
\definecolor{mycolor5}{rgb}{0.00000,0.72727,1.00000}%
\definecolor{mycolor6}{rgb}{0.933, 0.611, 0.054}%
\definecolor{mycolor7}{rgb}{0.90909,0.00000,1.00000}%
\definecolor{mycolor8}{rgb}{1.00000,0.00000,0.54545}%
\begin{tikzpicture}

\begin{axis}[%
width=\fwidth,
height=\fheight,
at={(0\fwidth,0\fheight)},
scale only axis,
xmin=0.005,
xmax=1,
xticklabel style={
	/pgf/number format/fixed,
	/pgf/number format/precision=3
},
yticklabel style={font=\tiny},
xtick={0.005,   0.1,   0.2,   0.3,   0.4,   0.5,   0.6,   0.7,   0.8,   0.9,     1},
ymin=0.5,
ymax=3.5,
ytick={1,2,3},
yticklabels={{\textsf{QDEIM}},{\textsf{KDEIM}},{\textsf{QR}}},
axis background/.style={fill=white}
]
\addplot [color=darkgray, line width=1.0pt,mark size=1.0pt, draw=none, mark=x, mark options={solid, darkgray}]
  table[row sep=crcr]{%
0.005	1\\
0.0150505050505051	1\\
0.0251010101010101	1\\
0.0351515151515151	1\\
0.0653030303030303	1\\
1	1\\
0.115555555555556	1\\
0.195959595959596	1\\
0.336666666666667	1\\
0.839191919191919	1\\
0.567828282828283	1\\
};
\addplot [color=blue!55, line width=0.5pt, draw=none, mark size=1.3pt, mark=square*, mark options={solid, blue!55}]
  table[row sep=crcr]{%
0.286414141414141	2\\
0.296464646464646	2\\
0.306515151515152	2\\
0.316565656565657	2\\
0.326616161616162	2\\
0.336666666666667	2\\
0.346717171717172	2\\
0.356767676767677	2\\
0.366818181818182	2\\
0.376868686868687	2\\
0.386919191919192	2\\
0.396969696969697	2\\
0.407020202020202	2\\
0.417070707070707	2\\
0.427121212121212	2\\
0.437171717171717	2\\
0.447222222222222	2\\
0.457272727272727	2\\
0.467323232323232	2\\
};
\addplot [color=mycolor4, line width=0.5pt, draw=none, mark size=1.3pt, mark=square*, mark options={solid, mycolor4}]
  table[row sep=crcr]{%
0.165808080808081	2\\
0.175858585858586	2\\
0.185909090909091	2\\
0.195959595959596	2\\
0.206010101010101	2\\
0.216060606060606	2\\
0.226111111111111	2\\
0.236161616161616	2\\
0.246212121212121	2\\
0.256262626262626	2\\
0.266313131313131	2\\
0.276363636363636	2\\
};
\addplot [color=mycolor8!55, line width=0.5pt, draw=none, mark size=1.3pt, mark=square*, mark options={solid, mycolor8!55}]
  table[row sep=crcr]{%
0.929646464646465	2\\
0.93969696969697	2\\
0.949747474747475	2\\
0.95979797979798	2\\
0.969848484848485	2\\
0.97989898989899	2\\
0.989949494949495	2\\
1	2\\
};
\addplot [color=mycolor6, line width=0.5pt, draw=none, mark size=1.3pt, mark=square*, mark options={solid, mycolor6}]
  table[row sep=crcr]{%
0.0150505050505051	2\\
};
\addplot [color=mycolor7, line width=0.5pt, draw=none, mark size=1.3pt, mark=square*, mark options={solid, mycolor7}]
  table[row sep=crcr]{%
0.005	2\\
};
\addplot [color=blue!20!mycolor2, line width=0.5pt, draw=none, mark size=1.3pt, mark=square*, mark options={solid, blue!20!mycolor2}]
  table[row sep=crcr]{%
0.708535353535354	2\\
0.718585858585859	2\\
0.728636363636364	2\\
0.738686868686869	2\\
0.748737373737374	2\\
0.758787878787879	2\\
0.768838383838384	2\\
0.778888888888889	2\\
0.788939393939394	2\\
0.798989898989899	2\\
0.809040404040404	2\\
0.819090909090909	2\\
0.829141414141414	2\\
0.839191919191919	2\\
0.849242424242424	2\\
0.859292929292929	2\\
0.869343434343434	2\\
0.879393939393939	2\\
0.889444444444444	2\\
0.899494949494949	2\\
0.909545454545455	2\\
0.91959595959596	2\\
};
\addplot [color=mycolor1, line width=0.5pt, draw=none, mark size=1.3pt, mark=square*, mark options={solid, mycolor1}]
  table[row sep=crcr]{%
0.0552525252525252	2\\
0.0653030303030303	2\\
0.0753535353535354	2\\
0.0854040404040404	2\\
};
\addplot [color=mycolor5, line width=0.5pt, draw=none, mark size=1.3pt, mark=square*, mark options={solid, mycolor5}]
  table[row sep=crcr]{%
0.0251010101010101	2\\
};
\addplot [color=mycolor2, line width=0.5pt, draw=none, mark size=1.3pt, mark=square*, mark options={solid, mycolor2}]
  table[row sep=crcr]{%
0.0351515151515151	2\\
0.0452020202020202	2\\
};
\addplot [color=mycolor3!60!white, line width=0.5pt, draw=none, mark size=1.3pt, mark=square*, mark options={solid, mycolor3!60!white}]
  table[row sep=crcr]{%
0.0954545454545455	2\\
0.105505050505051	2\\
0.115555555555556	2\\
0.125606060606061	2\\
0.135656565656566	2\\
0.145707070707071	2\\
0.155757575757576	2\\
};
\addplot [color=red!30, line width=0.5pt, draw=none, mark size=1.3pt, mark=square*, mark options={solid, red!30}]
  table[row sep=crcr]{%
0.477373737373737	2\\
0.487424242424242	2\\
0.497474747474748	2\\
0.507525252525252	2\\
0.517575757575758	2\\
0.527626262626263	2\\
0.537676767676768	2\\
0.547727272727273	2\\
0.557777777777778	2\\
0.567828282828283	2\\
0.577878787878788	2\\
0.587929292929293	2\\
0.597979797979798	2\\
0.608030303030303	2\\
0.618080808080808	2\\
0.628131313131313	2\\
0.638181818181818	2\\
0.648232323232323	2\\
0.658282828282828	2\\
0.668333333333333	2\\
0.678383838383838	2\\
0.688434343434343	2\\
0.698484848484848	2\\
};
\addplot [color=darkgray, line width=1.0pt, draw=none, mark size=1.0pt, mark=x, mark options={solid, darkgray}]
  table[row sep=crcr]{%
0.376868686868687	2\\
0.216060606060606	2\\
0.969848484848485	2\\
0.0150505050505051	2\\
0.005	2\\
0.809040404040404	2\\
0.0653030303030303	2\\
0.0251010101010101	2\\
0.0351515151515151	2\\
0.125606060606061	2\\
0.587929292929293	2\\
};
\addplot [color=darkgray, line width=1.0pt,mark size=1.0pt, draw=none, mark=x, mark options={solid, darkgray}]
  table[row sep=crcr]{%
0.005	3\\
0.748737373737374	3\\
0.105505050505051	3\\
0.0251010101010101	3\\
0.296464646464646	3\\
1	3\\
0.0150505050505051	3\\
0.0552525252525252	3\\
0.185909090909091	3\\
0.0351515151515151	3\\
0.507525252525252	3\\
0.899494949494949	3\\
};
\end{axis}
\end{tikzpicture}%
	\caption{\cref{alg:rbadapt2} for the Burgers' Equation with SVD, QR tolerance $\epsilon_{\text{SVD}}, \epsilon_{\text{QR}} = 1\cdot10^{-6}$. The crossmarks denote the parameters in the subsampled training set. For \textsf{KDEIM} each colour represents one cluster; the centroids of each of the clusters make up the subsampled training set.}
	\label{fig:burgtrngset_minus6}
\end{figure}
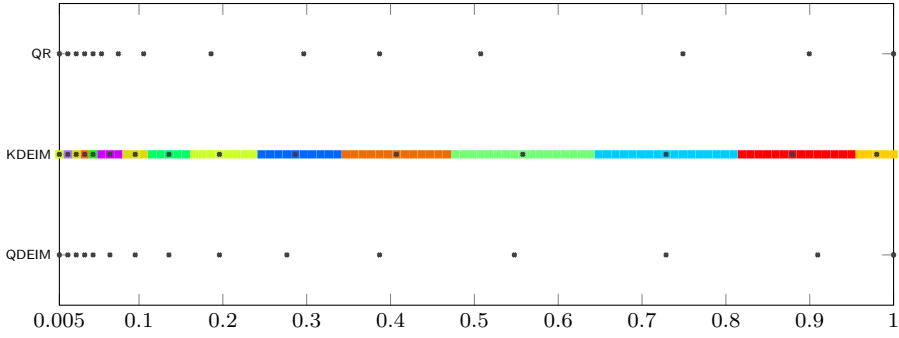
\begin{figure}[tp!]
	\centering
\definecolor{mycolor1}{rgb}{1.00000,0.80000,0.00000}%
\definecolor{mycolor2}{rgb}{0.80000,1.00000,0.00000}%
\definecolor{mycolor3}{rgb}{0.00000,1.00000,0.40000}%
\definecolor{mycolor4}{rgb}{0.00000,1.00000,0.80000}%
\definecolor{mycolor5}{rgb}{0.00000,0.80000,1.00000}%
\definecolor{mycolor6}{rgb}{0.40000,0.00000,1.00000}%
\definecolor{mycolor7}{rgb}{0.80000,0.00000,1.00000}%
\definecolor{mycolor8}{rgb}{1.00000,0.00000,0.80000}%
\definecolor{mycolor8}{rgb}{0.874, 0.850, 0.031}
\begin{tikzpicture}

\begin{axis}[%
width=\fwidth,
height=\fheight,
at={(0\fwidth,0\fheight)},
scale only axis,
xmin=0.005,
xmax=1,
xticklabel style={
	/pgf/number format/fixed,
	/pgf/number format/precision=3
},
yticklabel style={font=\tiny},
xtick={0.005,   0.1,   0.2,   0.3,   0.4,   0.5,   0.6,   0.7,   0.8,   0.9,     1},
ymin=0.5,
ymax=3.5,
ytick={1,2,3},
yticklabels={{\textsf{QDEIM}},{\textsf{KDEIM}},{\textsf{QR}}},
axis background/.style={fill=white}
]
\addplot [color=darkgray, line width=1.0pt,mark size=1.0pt,  draw=none, mark=x, mark options={solid, darkgray}]
  table[row sep=crcr]{%
0.005	1\\
0.0150505050505051	1\\
0.0251010101010101	1\\
0.0351515151515151	1\\
0.0452020202020202	1\\
0.0653030303030303	1\\
1	1\\
0.0954545454545455	1\\
0.135656565656566	1\\
0.195959595959596	1\\
0.276363636363636	1\\
0.909545454545455	1\\
0.386919191919192	1\\
0.547727272727273	1\\
0.728636363636364	1\\
};
\addplot [color=red, line width=0.5pt, draw=none, mark size=1.3pt, mark=square*, mark options={solid, red}]
  table[row sep=crcr]{%
0.819090909090909	2\\
0.829141414141414	2\\
0.839191919191919	2\\
0.849242424242424	2\\
0.859292929292929	2\\
0.869343434343434	2\\
0.879393939393939	2\\
0.889444444444444	2\\
0.899494949494949	2\\
0.909545454545455	2\\
0.91959595959596	2\\
0.929646464646465	2\\
0.93969696969697	2\\
0.949747474747475	2\\
};
\addplot [color=red!20!orange, line width=0.5pt, draw=none, mark size=1.3pt, mark=square*, mark options={solid, red!20!orange}]
  table[row sep=crcr]{%
0.0351515151515151	2\\
};
\addplot [color=mycolor1, line width=0.5pt, draw=none, mark size=1.3pt, mark=square*, mark options={solid, mycolor1}]
  table[row sep=crcr]{%
0.95979797979798	2\\
0.969848484848485	2\\
0.97989898989899	2\\
0.989949494949495	2\\
1	2\\
};
\addplot [color=mycolor2, line width=0.5pt, draw=none, mark size=1.3pt, mark=square*, mark options={solid, mycolor2}]
  table[row sep=crcr]{%
0.005	2\\
};
\addplot [color=mycolor8, line width=0.5pt, draw=none, mark size=1.3pt, mark=square*, mark options={solid, mycolor8}]
  table[row sep=crcr]{%
0.0251010101010101	2\\
};
\addplot [color=green, line width=0.5pt, draw=none, mark size=1.3pt, mark=square*, mark options={solid, green}]
  table[row sep=crcr]{%
0.0452020202020202	2\\
};
\addplot [color=mycolor3, line width=0.5pt, draw=none, mark size=1.3pt, mark=square*, mark options={solid, mycolor3}]
  table[row sep=crcr]{%
0.115555555555556	2\\
0.125606060606061	2\\
0.135656565656566	2\\
0.145707070707071	2\\
0.155757575757576	2\\
};
\addplot [color=mycolor4!20!yellow, line width=0.5pt, draw=none, mark size=1.3pt, mark=square*, mark options={solid, mycolor4!20!yellow}]
  table[row sep=crcr]{%
0.165808080808081	2\\
0.175858585858586	2\\
0.185909090909091	2\\
0.195959595959596	2\\
0.206010101010101	2\\
0.216060606060606	2\\
0.226111111111111	2\\
0.236161616161616	2\\
};
\addplot [color=mycolor5, line width=0.5pt, draw=none, mark size=1.3pt, mark=square*, mark options={solid, mycolor5}]
  table[row sep=crcr]{%
0.648232323232323	2\\
0.658282828282828	2\\
0.668333333333333	2\\
0.678383838383838	2\\
0.688434343434343	2\\
0.698484848484848	2\\
0.708535353535354	2\\
0.718585858585859	2\\
0.728636363636364	2\\
0.738686868686869	2\\
0.748737373737374	2\\
0.758787878787879	2\\
0.768838383838384	2\\
0.778888888888889	2\\
0.788939393939394	2\\
0.798989898989899	2\\
0.809040404040404	2\\
};
\addplot [color=blue!50!mycolor5, line width=0.5pt, draw=none, mark size=1.3pt, mark=square*, mark options={solid, blue!50!mycolor5}]
  table[row sep=crcr]{%
0.246212121212121	2\\
0.256262626262626	2\\
0.266313131313131	2\\
0.276363636363636	2\\
0.286414141414141	2\\
0.296464646464646	2\\
0.306515151515152	2\\
0.316565656565657	2\\
0.326616161616162	2\\
0.336666666666667	2\\
};
\addplot [color=green!55, line width=0.5pt, draw=none, mark size=1.3pt, mark=square*, mark options={solid, green!55}]
  table[row sep=crcr]{%
0.477373737373737	2\\
0.487424242424242	2\\
0.497474747474748	2\\
0.507525252525252	2\\
0.517575757575758	2\\
0.527626262626263	2\\
0.537676767676768	2\\
0.547727272727273	2\\
0.557777777777778	2\\
0.567828282828283	2\\
0.577878787878788	2\\
0.587929292929293	2\\
0.597979797979798	2\\
0.608030303030303	2\\
0.618080808080808	2\\
0.628131313131313	2\\
0.638181818181818	2\\
};
\addplot [color=mycolor6!55, line width=0.5pt, draw=none, mark size=1.3pt, mark=square*, mark options={solid, mycolor6!55}]
  table[row sep=crcr]{%
0.0150505050505051	2\\
};
\addplot [color=mycolor7, line width=0.5pt, draw=none, mark size=1.3pt, mark=square*, mark options={solid, mycolor7}]
  table[row sep=crcr]{%
0.0552525252525252	2\\
0.0653030303030303	2\\
0.0753535353535354	2\\
};
\addplot [color=mycolor8, line width=0.5pt, draw=none, mark size=1.3pt, mark=square*, mark options={solid, mycolor8}]
  table[row sep=crcr]{%
0.0854040404040404	2\\
0.0954545454545455	2\\
0.105505050505051	2\\
};
\addplot [color=red!50!mycolor8, line width=0.5pt, draw=none, mark size=1.3pt, mark=square*, mark options={solid, red!50!mycolor8}]
  table[row sep=crcr]{%
0.346717171717172	2\\
0.356767676767677	2\\
0.366818181818182	2\\
0.376868686868687	2\\
0.386919191919192	2\\
0.396969696969697	2\\
0.407020202020202	2\\
0.417070707070707	2\\
0.427121212121212	2\\
0.437171717171717	2\\
0.447222222222222	2\\
0.457272727272727	2\\
0.467323232323232	2\\
};
\addplot [color=darkgray, line width=1.0pt,mark size=1.0pt, draw=none, mark size=1.0pt, mark=x, mark options={solid, darkgray}]
  table[row sep=crcr]{%
0.879393939393939	2\\
0.0351515151515151	2\\
0.97989898989899	2\\
0.005	2\\
0.0251010101010101	2\\
0.0452020202020202	2\\
0.135656565656566	2\\
0.195959595959596	2\\
0.728636363636364	2\\
0.286414141414141	2\\
0.557777777777778	2\\
0.0150505050505051	2\\
0.0653030303030303	2\\
0.0954545454545455	2\\
0.407020202020202	2\\
};
\addplot [color=darkgray, line width=1.0pt,mark size=1.0pt, draw=none, mark=x, mark options={solid, darkgray}]
  table[row sep=crcr]{%
0.005	3\\
0.748737373737374	3\\
0.105505050505051	3\\
0.0251010101010101	3\\
0.296464646464646	3\\
1	3\\
0.0150505050505051	3\\
0.0552525252525252	3\\
0.185909090909091	3\\
0.0351515151515151	3\\
0.507525252525252	3\\
0.899494949494949	3\\
0.0753535353535354	3\\
0.386919191919192	3\\
0.0452020202020202	3\\
};
\end{axis}
\end{tikzpicture}%
	\caption{\cref{alg:rbadapt2} for the Burgers' Equation with SVD, QR tolerance $\epsilon_{\text{SVD}}, \epsilon_{\text{QR}} = 1\cdot10^{-8}$. The crossmarks denote the parameters in the subsampled training set. For \textsf{KDEIM} each colour represents one cluster; the centroids of each of the clusters make up the subsampled training set.}
	\label{fig:burgtrngset_minus8}
\end{figure}
\setlength\fheight{6cm}
\setlength\fwidth{6cm}	
\begin{figure}[t]
	\centering
	\subfloat[\textsf{QDEIM}.]{\label{fig:err_plot_alg41_qdeim}  \begin{tikzpicture}
    \begin{axis}[
    	view   = {0}{90},
        axis on top,% ----
        width=\fwidth,
        height=\fheight,
        scale only axis,
        enlargelimits=false,
        xmin=0,
        xmax=2,
        ymin=0.005,
        ymax=1,
        xlabel={Time (s)},
        ylabel={Parameter $\bm{\mu}$},
        axis on top, 
        axis equal image,
		point meta min=0,
		point meta max=4.5463e-08,
		colormap/viridis,
		colorbar
        ]
      \addplot graphics[xmin=0,ymin=0.005,xmax=2,ymax=1] {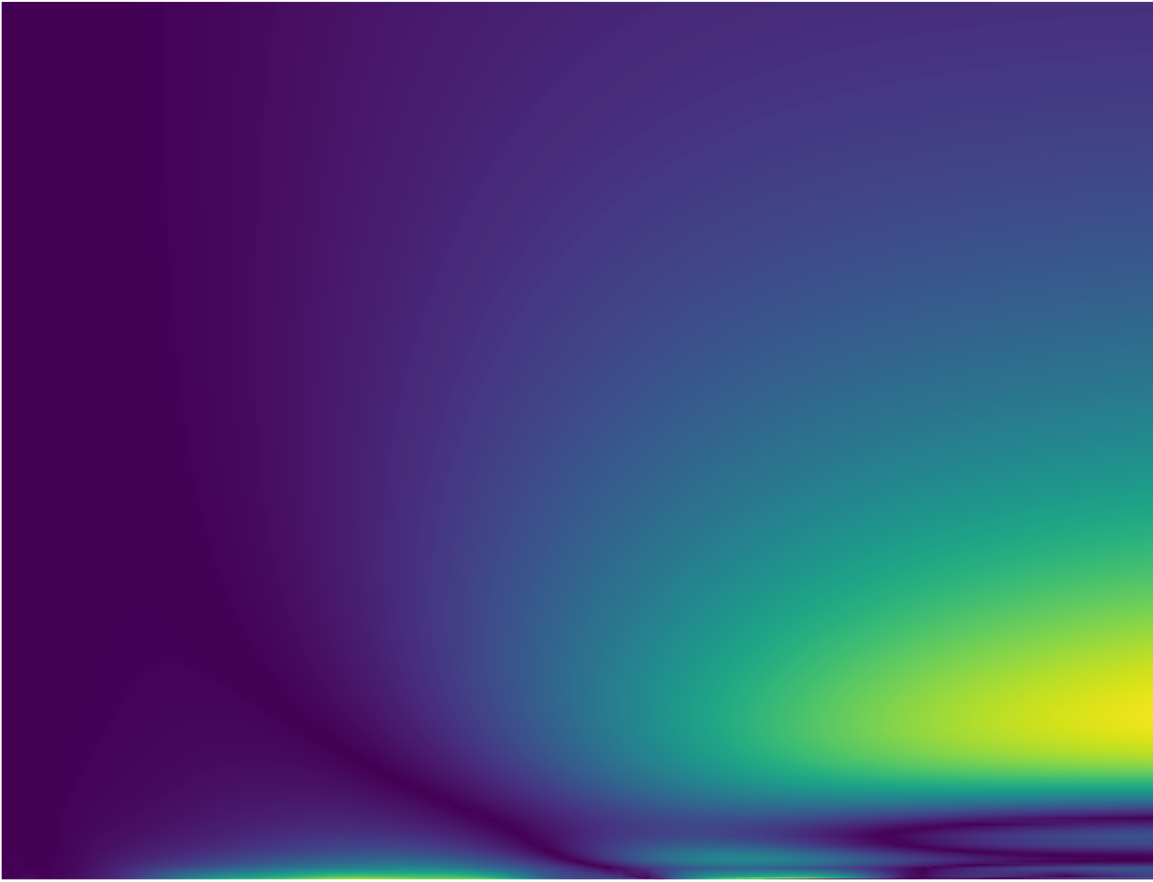};
    \end{axis}
  \end{tikzpicture}}\hfill
	\subfloat[\textsf{KDEIM}.]{\label{fig:err_plot_alg41_kdeim}  \begin{tikzpicture}
    \begin{axis}[
    	view   = {0}{90},
        axis on top,% ----
        width=\fwidth,
        height=\fheight,
        scale only axis,
        enlargelimits=false,
        xmin=0,
        xmax=2,
        ymin=0.005,
        ymax=1,
        xlabel={Time (s)},
		ylabel={Parameter $\bm{\mu}$},        
        axis on top, 
        axis equal image,
		point meta min=0,
		point meta max=6.8796e-08,
		colormap/viridis,
		colorbar
        ]
      \addplot graphics[xmin=0,ymin=0.005,xmax=2,ymax=1] {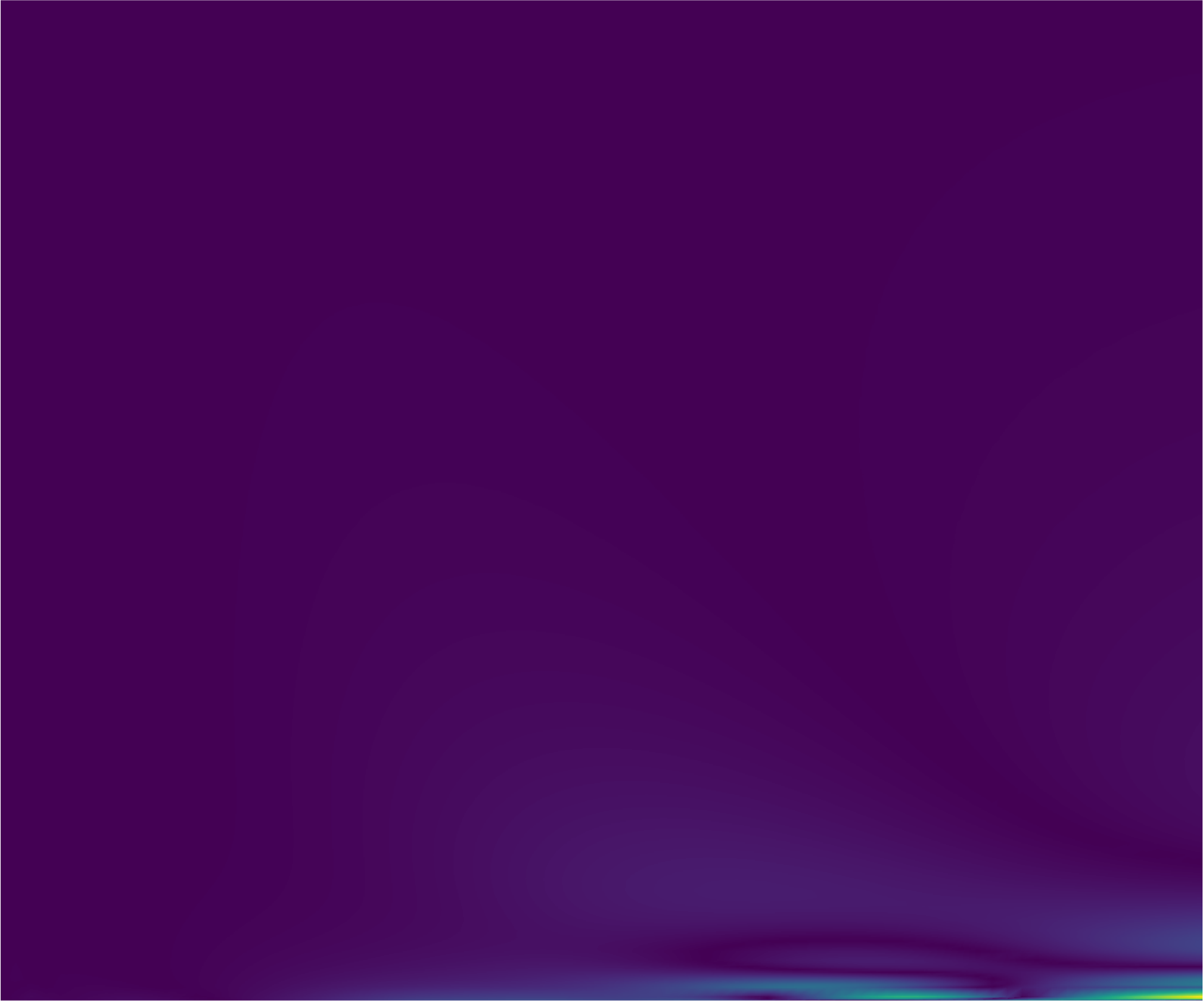};
    \end{axis}
  \end{tikzpicture}}
	\caption{Error plot for \cref{alg:rbadapt1} with tolerance $\epsilon_{\text{SVD}} = 10^{-8}$ applied to the Burgers' equation. The error between the true and reduced outputs $\| \bx{y}\left(\bx{x}(t^{k},\bm{\mu})\right) - \widetilde{\bx{y}}\left(\bx{z}(t^{k},\bm{\mu})\right) \|$ is plotted over the duration of the simulation for all parameters in the test set.}
	\label{fig:err_plot_alg41}
\end{figure}
\subsection{Thermal Block}
\begin{figure}[t!]
	\centering
	\includegraphics[scale = 0.5]{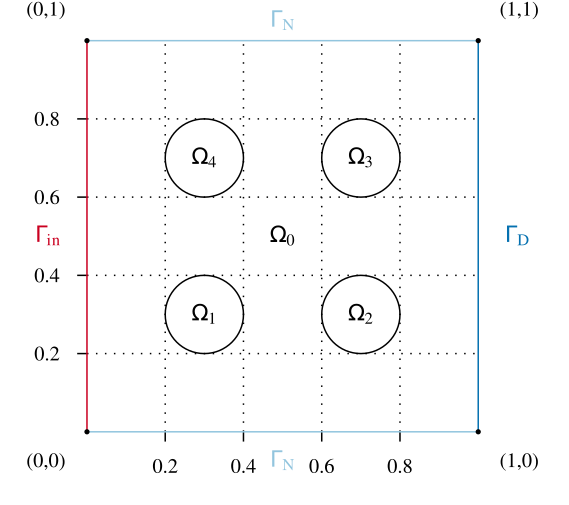}
	\caption{Thermal Block Example: Spatial domain and boundaries.}
	\label{fig:cookie}
\end{figure}

The second example is a benchmark model of the time-dependent heat transfer in a thermal block. The governing PDE is given by
\begin{equation}
	\begin{aligned}
		\frac{\partial \theta(\bx{w}, t, \bm{\mu})}{\partial t} + \nabla \cdot (-\gamma(\bx{w}, \bm{\mu})\nabla\theta(\bx{w}, t, \bm{\mu})) &= 0, \qquad t \in [0,\,T].
	\end{aligned}
\end{equation}
The domain $\Omega := (0 , 1) \times (0 , 1) \in \R^{2}$ is partitioned into five regions --- $\Omega = \Omega_{0} \cup \Omega_{1} \cup \Omega_{2} \cup \Omega_{3} \cup \Omega_{4}$ as shown in \cref{fig:cookie}. The left boundary of the domain ($\Gamma_{\text{in}}$) is associated with an input heat flux of unit magnitude, the top and bottom boundaries ($\Gamma_{\text{N}}$) are associated with a Neumann boundary condition with zero flux and finally the right boundary ($\Gamma_{\text{D}}$) is fixed at zero. The state variable is the temperature $\theta(\bx{w}, t)$ at a given spatial location $\bx{w} \in \Omega$, for a given time $t$. The initial condition is $\theta(\bx{w},0) = 0$. The output is the average temperature measured at $\Omega_{2}$. The problem is parametrized by the heat conductivity $\gamma$ in the subdomains ($\Omega_{0}, \Omega_{1}, \Omega_{2}, \Omega_{3}$ and $\Omega_{4}$); $\gamma(\bx{w}, \bm{\mu}) = 1$ when $\bx{w} \in \Omega_{0}$ and $\gamma(\bx{w}, \bm{\mu}) = \kappa_{i}$ whenever $\bx{w} \in \Omega_{i}$, $i = 1,2,3,4$. 
We define the parameter vector $\bm{\mu} = [\kappa_{1}, \kappa_{2}, \kappa_{3}, \kappa_{4}]$. The governing PDE is discretized in space using linear finite elements with respect to a simplicial triangulation of the domain $\Omega$ obtained via the software \texttt{gmsh} \cite{gmsh09}. It is further discretized in time using the implicit Euler scheme for a time ranging from $t \in [0 \,, 1]$, with step size $\Delta t = 0.01$. The spatially discretized system has dimension $n = 7488$. For more details on the model and the spatial discretization, the reader is referred to \cite{morRavS20}. The discretized heat equation can be written in the form of \cref{eq:odess}. Since the problem is linear, we have $\bx{f} \equiv 0$. For the numerical results, the parameter $\bm{\mu}$ is sampled from the domain $\mathcal{P} := [{1\cdot10^{-5}} \,, {1\cdot10^{-2}}] \times [{1\cdot10^{-5}} \,, {1\cdot10^{-2}}] \times [{1\cdot10^{-4}} \,, 1] \times [{1\cdot10^{-1}} \,, 1]$. For purposes of illustration, we consider the three parameter version of the thermal block problem by fixing $\kappa_{4}$ to its mean value, i.e., $\kappa_{4} = 0.5$. The training set $\Xi_{\text{train}}$ consists of a tensor grid of $N_{\text{train}} = 6^{3} = 216$ parameters, with $6$ parameters sampled for each $\kappa_{i}, \,\, i = 1,\, 2,\, 3$. The test parameter set consists of 100 parameters, randomly sampled from $\mathcal{P}$. The tolerance for the greedy algorithm is set to be $\texttt{tol} = 1\cdot10^{-3}$.

\subsubsection{Greedy Algorithm with Fixed Training Set}
\setlength\fheight{4cm}
\setlength\fwidth{6cm}		
\begin{figure}[t]
\centering
\input{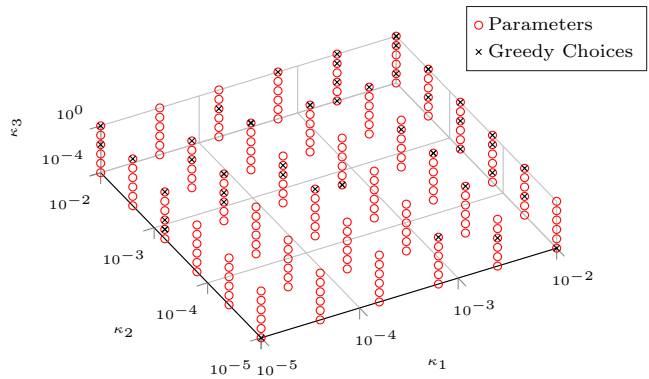}
\caption{Thermal Block: Fine training set with $216$ parameters and the $44$ greedy parameters picked by \cref{alg:rb}.}
\label{fig:thermal_fixed}
\end{figure}
Applying \cref{alg:rb} with a fixed training set to the thermal block example results in a ROM of dimension $r_{\text{POD}} = 74$, taking $53$ iterations to converge in $t_{\text{greedy}} = 694.73$ seconds. The maximum error over the test set is $\epsilon_{t}^{\text{max}} = 9.78 \cdot 10^{-4}$. In \cref{fig:thermal_fixed}, the training set $\Xi_{\textnormal{train}}$ and the greedy parameters identified by \cref{alg:rb} are shown. Of the $216$ parameters in the training set, only $44$ are chosen. The greedy parameters have a larger concentration at and around $(0.01, 0.01, 0.0001)$, the upper right corner of the figure. In fact, the regions around the vicinity of the upper and right wall of the grid posses many greedy samples near them. 

\subsubsection{Greedy Algorithm Schemes 1 and 2}
Similar to the Burgers' equation, we now apply the proposed training set subsampling schemes to the thermal block example. We shall also illustrate the advantages of using oversampling. For both \cref{alg:rbadapt1,alg:rbadapt2}, we consider $\epsilon_{\text{QR}}, \epsilon_{\text{SVD}} = 1\cdot10^{-10}$ and a coarse tolerance $\texttt{tol}^{c} = 1$. The approximation to $\mathbf{\mathfrak{Y}}$ is obtained by taking snapshots at every time step of the implicit Euler scheme. The results are summarized in \cref{tab:thermal_v1} and \cref{tab:thermal_v2} for \cref{alg:rbadapt1,alg:rbadapt2}, respectively. The first scheme does not lead to a successful ROM for both \textsf{QDEIM} and \textsf{KDEIM} whereas using the pivoted QR decomposition on $\widetilde{\mathbf{\mathfrak{Y}}}^{\textsf{T}}$ to identify the subsampled training set produces a successful ROM. For the second scheme of the proposed algorithm, both \textsf{QDEIM} and \textsf{KDEIM}  result in a subsampled training set of cardinality $N_{c} = 19$ while the pivoted QR approach gives $N_{c} = 20$. However, both QR and \textsf{QDEIM} are unsuccessful in meeting the required ROM tolerance for the test set. On the other hand, \textsf{KDEIM} results in a ROM satisfying the tolerance, taking a significantly smaller number of iterations ($40$) to converge. The results seem to indicate that the subsampling approach is not entirely able to capture the full range of features over the training set. This is mainly due to the smaller number of parameters $N_{c} = 19$, that the algorithm results in. Recall that for the standard greedy approach, $44$ unique greedy parameters were determined. However, it is also to be noted that the performance of the ROMs resulting from either scheme on the test set is not bad. The maximum error is only slightly higher than the desired tolerance of $\texttt{tol}  = 1\cdot10^{-3}$. 

\begin{figure}[t!]
\centering
\subfloat[Subsampled parameters using \textsf{Gappy-POD Eigenvector}.]{\label{fig:oversample_gpode}\input{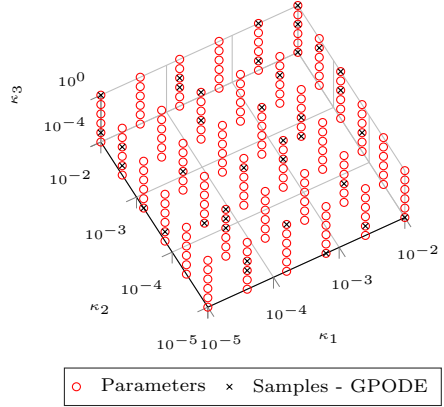}}\hfill	
\subfloat[Subsampled parameters using \textsf{Gappy-POD Clustering}.]{\label{fig:oversample_gpodc}\input{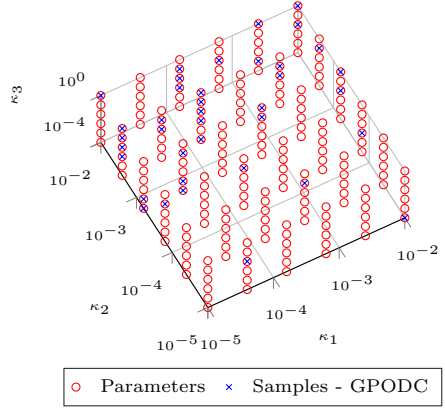}} 
\caption{Subsampling strategy using \textsf{Gappy-POD} with oversampling for the Thermal Block.}
\label{fig:oversample}
\end{figure}

Next, we perform oversampling to identify more parameters. For the standard \textsf{DEIM} approach, the number of interpolation points is equal to the rank $\ell$ of the truncated left singular vectors of the snapshots matrix. For oversampling, we set the number of interpolation points to be $m = 2 \ell$ and test the approaches based on maximizing the smallest singular value (\textsf{Gappy-POD Eigenvector}) and the approach based on clustering (\textsf{Gappy-POD Clustering}), both originally proposed in \cite{morPehDG18}. The results are given in \cref{tab:thermal_oversample}. We see that both oversampling approaches result in ROMs that are validated to be accurate over the test set. The \textsf{Gappy-POD Clustering} method results in the smallest test error among the two and takes $40$ iterations to converge. Notice that, compared to the previous two approaches based on \textsf{QDEIM} and \textsf{KDEIM}, the oversampling approach requires more time. This is not surprising, since a larger number of parameters is included in the coarse training set at Stage 2 of \cref{alg:rbadapt1,alg:rbadapt2}. The speedup of the \textsf{Gappy-POD Eigenvector} variant is $3.9$, while a speedup of $4.6$ is achieved by the \textsf{Gappy-POD Clustering} variant. We show the subsampled training sets of both the approaches in \cref{fig:oversample}. In particular, parameter samples anticipated by the \textsf{Gappy-POD Clustering} variant bear a close resemblance to the greedy parameter distribution in \cref{fig:thermal_fixed}. In \cref{fig:errplot_gpod}, we plot the mean error over time for each parameter in the test set. It is evident that both the proposed oversampling strategies, \textsf{Gappy-POD Eigenvector} (\cref{fig:errplot_gpode}) and \textsf{Gappy-POD Clustering} (\cref{fig:errplot_gpodc})  have been successful in keeping the error below the desired tolerance, uniformly for all the parameters in the test set. Lastly, for the thermal block example, the total number of iterations to converge for \cref{alg:rb}, \cref{alg:rbadapt1,alg:rbadapt2} is not equal and varies based on the particular subsampling strategy used (see \cref{tab:thermal_v1,tab:thermal_v2,tab:thermal_oversample}). Our initial assumptions in \cref{subsec:compl_analysis} that $N_{I} \approx N_{II,1} + N_{II,2}$, $N_{I} \approx N_{III,1} + N_{III,2}$ do not hold while we see that $N_{I} \gtrapprox N_{II,1} + N_{II,2}$, $N_{I} \gtrapprox N_{III,1} + N_{III,2}$ is always the case.
 
% Thermal Block Table 1
\begin{table}[t!]
\caption{Thermal Block Results for \cref{alg:rbadapt1} for \textsf{QDEIM, KDEIM} and QR.}
\label{tab:thermal_v1}
\centering
\small
\begin{tabular}{|c|c|c|c|c|}
\hline
\multirow{3}{*}{Method}    & \multirow{3}{*}{Fixed} & \multicolumn{3}{c|}{Adapted}                                                                                    \\ \cline{3-5} 
                           &                                     & \multicolumn{3}{c|}{$\epsilon_{\text{SVD}}, \epsilon_{\text{QR}} = 1\cdot10^{-10}$}\\ \cline{3-5} 
                           &                                     & \small{\textsf{QDEIM}}              & \small{\textsf{KDEIM}} & \small{QR}\\ \hline
$N_{c}$                    	& 216                                 & 19                 & 19 	& 20			\\ \hline
$\epsilon_{t}^{\text{max}}$ & $9.78 \cdot 10^{-4}$                & $1.10\cdot10^{-3}$ & $2.10\cdot10^{-3}$ & $6.45\cdot10^{-4}$\\ \hline
$r_{\text{POD}}$          & 74                             	  & 64                 & 64       & 66     		\\ \hline
Iterations (\texttt{iter})                       & 53                                  & 42                 & 43      & 44    		\\ \hline
Offline time (s)            & 694.73                              & 121.36             & 123.27 	  		& 126.35 \\ \hline
Speedup                     & -                                   & 5.7                & 5.6          	& 5.5	\\ \hline
\end{tabular}
\end{table}
% Thermal Block Table 2
\begin{table}[t!]
\caption{Thermal Block Results for \cref{alg:rbadapt2} for \textsf{QDEIM, KDEIM} and QR.}
\label{tab:thermal_v2}
\centering
\small
\begin{tabular}{|c|c|c|c|c|}
\hline
\multirow{3}{*}{Method}    & \multirow{3}{*}{Fixed} & \multicolumn{3}{c|}{Adapted}                                                                                    \\ \cline{3-5} 
                           &                                     & \multicolumn{3}{c|}{$\epsilon_{\text{SVD}}, \epsilon_{\text{QR}} = 1\cdot10^{-10}$}\\ \cline{3-5} 
                           &                                     & \small{\textsf{QDEIM}}              & \small{\textsf{KDEIM}} & \small{QR}\\ \hline
$N_{\text{train}}$                    	& 216                                 & 19                 & 19 			& 20	\\ \hline
$\epsilon_{t}^{\text{max}}$ & $9.78 \cdot 10^{-4}$                & $1.10\cdot10^{-3}$ & $9.36\cdot10^{-4}$ & $1.60\cdot10^{-3}$\\ \hline
$r_{\text{POD}}$          & 74                             	  & 65                 & 62            & 57		\\ \hline
Iterations (\texttt{iter})                       & 53                                  & 45                 & 40   & 35       		\\ \hline
Offline time (s)            & 694.73                              & 121.58             & 110.13 	  & 73.38		\\ \hline
Speedup                     & -                                   & 5.7                & 6.3        & 9.5  		\\ \hline
\end{tabular}
\end{table}
% Thermal Block Oversampling
\begin{table}[t!]
\caption{Thermal Block Results for \cref{alg:rbadapt1} with oversampling.}
\label{tab:thermal_oversample}
\centering
\small
\begin{tabular}{|c|c|c|}
\hline
\multirow{3}{*}{Method}   & \multicolumn{2}{c|}{Oversampling}                                                                                    \\ \cline{2-3} 
                          &  \multicolumn{2}{c|}{$m = 2 \ell$}\\ \cline{2-3} 
                          & \small{\textsf{Gappy-POD Eigenvector}}              & \small{\textsf{Gappy-POD Clustering}}\\ \hline
$N_{c}$                    	&  38                 & 38 					\\ \hline
$\epsilon_{t}^{\text{max}}$ &  $9.91\cdot10^{-4}$ & $7.96\cdot10^{-4}$ 	\\ \hline
$r_{\text{POD}}$            &  70                 & 62            		\\ \hline
Iterations (\texttt{iter})                       &  47                 & 40          		\\ \hline
Offline time (s)            &  177.51             & 151.39 	  			\\ \hline
Speedup                     &  3.9                & 4.6          		\\ \hline
\end{tabular}
\end{table}
\setlength\fheight{4cm}
\setlength\fwidth{4cm}		
\setlength\fheight{4cm}
\setlength\fwidth{4cm}		
\begin{figure}[t!]
	\centering
	\subfloat[\textsf{Gappy-POD Eigenvector}.]{\label{fig:errplot_gpode}\input{figures/Ex2_Thermal_ErrPlot_v2_1_GPODE.tex}}\hfill	
	\subfloat[\textsf{Gappy-POD Clustering}.]{\label{fig:errplot_gpodc}\input{figures/Ex2_Thermal_ErrPlot_v2_1_GPODC.tex}} 
	\caption{Error plot for \cref{alg:rbadapt1} with coarse tolerance $\texttt{tol}^{c} = 1$ and subsampling based on \textsf{Gappy-POD Eigenvector} and \textsf{Gappy-POD Clustering} applied to the thermal block example. The mean error over time between the true and reduced outputs - $(1/K+1) \sum_{i = 0}^{K } \| \bx{y}\left(\bx{x}(t^{k},\bm{\mu})\right) - \widetilde{\bx{y}}\left(\bx{z}(t^{k},\bm{\mu})\right) \|$ - is plotted for all parameters in the test set.}
	\label{fig:errplot_gpod}
\end{figure}

\section{Conclusions}
\label{sec:conclusion}
We presented an efficient method to subsample the training set in the offline stage of the RBM. The proposed two-stage strategy is goal-oriented. It uses the pivoted QR decomposition, the \textsf{DEIM} algorithm or its variants to approximate the parameter-to-output map for time-dependent problems,  taking advantage of the information from the pivots or the interpolation points to generate the subsampled training set. Different strategies to identify the interpolation points based on variants of the \textsf{DEIM} algorithm were discussed. The strategy of retaining as many parameters as the number of \textsf{DEIM} interpolation points is simple and leads to considerable speedup. However, there may be occasions where it is not robust, as demonstrated in the thermal block example. For such scenarios, we propose a principled oversampling strategy based on the \textsf{Gappy-POD} approach, to add additional parameters to the subsampled training set. This approach, while slightly more expensive, is robust and leads to reliable ROMs. The different subsampling strategies were tested on two numerical examples and were shown to yield, on average, a speedup of up to $5$ for the offline stage of the RBM, without compromising the quality of the generated ROMs.

There exist several promising possibilities for further applications and improvements. One possible idea is in relation to the \textsf{KDEIM}-based subsampling strategy. It is worthwhile to consider the development of localized RBMs for each cluster identified by the algorithm. Such an approach could help overcome existing challenges related to training set adaptation methods that partition the parameter domain based on binary trees \cite{morEftKP11,morHaaDO11}. Another promising idea is related to data assimilation within the model-based RBM framework. Assuming output data is available by means of sensors or other measurements, it should be possible to consider a QR decomposition of this parametric data matrix to identify a good initial training set. Such an approach is capable of incorporating data in a natural fashion, to aid the  development of problem-tailored ROMs.
\section*{Code Availability}
The companion code and data used to compute the results are available at
\begin{center}
	\url{https://doi.org/10.5281/zenodo.4593144}
\end{center}
under the MIT License.
%\begin{acknowledgements}
%\end{acknowledgements}

\bibliographystyle{spmpsci}      % mathematics and physical sciences
\bibliography{references}

\begin{thebibliography}{10}
\providecommand{\url}[1]{{#1}}
\providecommand{\urlprefix}{URL }
\expandafter\ifx\csname urlstyle\endcsname\relax
  \providecommand{\doi}[1]{DOI~\discretionary{}{}{}#1}\else
  \providecommand{\doi}{DOI~\discretionary{}{}{}\begingroup
  \urlstyle{rm}\Url}\fi

\bibitem{Aan09}
Aanonsen, T.: Empirical interpolation with application to reduced basis
  approximations.
\newblock Master's thesis, Master Thesis, Norwegian University of Science and
  Technology (NTNU), Trondheim, Norway (2009).
\newblock \urlprefix\url{http://hdl.handle.net/11250/258487}

\bibitem{morAntCF18}
Antil, H., Chen, D., Field, S.: A note on qr-based model reduction: Algorithm,
  software, and gravitational wave applications.
\newblock Comput. Sci. Eng. \textbf{20}(4), 10--25 (2018).
\newblock \doi{10.1109/MCSE.2018.042781323}

\bibitem{morAnt05}
Antoulas, A.C.: Approximation of Large-Scale Dynamical Systems, \emph{Adv. Des.
  Control}, vol.~6.
\newblock {SIAM} Publications, Philadelphia, PA (2005).
\newblock \doi{10.1137/1.9780898718713}

\bibitem{AntBG20}
Antoulas, A.C., Beattie, C.A., Gugercin, S.: Interpolatory Methods for Model
  Reduction.
\newblock Society for Industrial and Applied Mathematics, Philadelphia, PA
  (2020).
\newblock \doi{10.1137/1.9781611976083}

\bibitem{morBarMNetal04}
Barrault, M., Maday, Y., Nguyen, N.C., Patera, A.T.: An `empirical
  interpolation' method: application to efficient reduced-basis discretization
  of partial differential equations.
\newblock C.R. Acad. Sci. Paris \textbf{339}(9), 667--672 (2004).
\newblock \doi{10.1016/j.crma.2004.08.006}

\bibitem{morBenEEetal18}
Benaceur, A., Ehrlacher, V., Ern, A., Meunier, S.: A progressive reduced
  basis/empirical interpolation method for nonlinear parabolic problems.
\newblock {SIAM} J. Sci. Comput. \textbf{40}(5), A2930--A2955 (2018).
\newblock \doi{10.1137/17M1149638}

\bibitem{morBenCOetal17}
Benner, P., Cohen, A., Ohlberger, M., Willcox, K. (eds.): Model Reduction and
  Approximation: Theory and Algorithms.
\newblock Computational Science \& Engineering. SIAM Publications,
  Philadelphia, PA (2017).
\newblock \doi{10.1137/1.9781611974829}

\bibitem{morBenetal21}
Benner, P., Grivet-Talocia, S., Quarteroni, A., Rozza, G., Schilder
  W.~Silveira, L.M. (eds.): {M}odel {O}rder {R}eduction. Volume 3:
  Applications.
\newblock De Gruyter (2021).
\newblock \doi{10.1515/9783110499001}

\bibitem{morBenGW15}
Benner, P., Gugercin, S., Willcox, K.: A survey of model reduction methods for
  parametric systems.
\newblock SIAM Review \textbf{57}(4), 483--531 (2015).
\newblock \doi{10.1137/130932715}

\bibitem{BouMD09}
Boutsidis, C., Mahoney, M.W., Drineas, P.: An improved approximation algorithm
  for the column subset selection problem.
\newblock In: Proceedings of the {T}wentieth {A}nnual {ACM}-{SIAM} {S}ymposium
  on {D}iscrete {A}lgorithms, pp. 968--977. SIAM, Philadelphia, PA (2009)

\bibitem{BroBP10}
Broadbent, M.E., Brown, M., Penner, K.: Subset selection algorithms: Randomized
  vs. deterministic.
\newblock SIAM Undergrad. Res. Online \textbf{3}, 50--71 (2010).
\newblock \doi{https://dx.doi.org/10.1137/09S010435}

\bibitem{morCaretal11}
Carlberg, K., Amsallem, D., Avery, P., Zahr, M., Farhat, C.: The {GNAT}
  nonlinear model reduction method and its application to fluid dynamics
  problems.
\newblock 6th AIAA Theoretical Fluid Mechanics Conference, Honolulu pp. 1--24
  (2011).
\newblock \doi{10.2514/6.2011-3112}

\bibitem{CivM12}
\c{C}ivril, A., Magdon-Ismail, M.: Column subset selection via sparse
  approximation of {SVD}.
\newblock Theoret. Comput. Sci. \textbf{421}, 1--14 (2012).
\newblock \doi{10.1016/j.tcs.2011.11.019}

\bibitem{Cha87}
Chan, T.F.: Rank revealing {$QR$} factorizations.
\newblock Linear Algebra Appl. \textbf{88/89}, 67--82 (1987).
\newblock \doi{10.1016/0024-3795(87)90103-0}

\bibitem{morChaS10}
Chaturantabut, S., Sorensen, D.C.: Nonlinear model reduction via discrete
  empirical interpolation.
\newblock {SIAM} J. Sci. Comput. \textbf{32}(5), 2737--2764 (2010).
\newblock \doi{10.1137/090766498}

\bibitem{morCheFB20}
Chellappa, S., Feng, L., Benner, P.: Adaptive basis construction and improved
  error estimation for parametric nonlinear dynamical systems.
\newblock Internat. J. Numer. Methods Engrg. \textbf{121}(23), 5320--5349
  (2020).
\newblock \doi{10.1002/nme.6462}

\bibitem{morCheFB19}
Chellappa, S., Feng, L., Benner, P.: An adaptive sampling approach for the
  reduced basis method.
\newblock In: Realization and Model Reduction of Dynamical Systems - A
  Festschrift in Honor of the 70th Birthday of Thanos Antoulas. Springer
  (2021).
\newblock \urlprefix\url{https://arxiv.org/abs/1910.00298}.
\newblock Accepted March 2020

\bibitem{morCheG19}
Chen, P., Ghattas, O.: Hessian-based sampling for high-dimensional model
  reduction.
\newblock Int. J. Uncertain. Quantif. \textbf{9}(2), 103--121 (2019).
\newblock \doi{10.1615/Int.J.UncertaintyQuantification.2019028753}

\bibitem{morCheQR17}
Chen, P., Quarteroni, A., Rozza, G.: Reduced basis methods for uncertainty
  quantification.
\newblock SIAM/ASA J. Uncertain. Quantif. \textbf{5}(1), 813--869 (2017).
\newblock \doi{10.1137/151004550}

\bibitem{morConDW14}
Constantine, P.G., Dow, E., Wang, Q.: Active subspace methods in theory and
  practice: applications to kriging surfaces.
\newblock SIAM J. Sci. Comput. \textbf{36}(4), A1500--A1524 (2014).
\newblock \doi{10.1137/130916138}

\bibitem{morDrmG16}
Drma\v{c}, Z., Gugercin, S.: A new selection operator for the discrete
  empirical interpolation method---improved a priori error bound and
  extensions.
\newblock {SIAM} J. Sci. Comput. \textbf{38}(2), A631--A648 (2016).
\newblock \doi{10.1137/15M1019271}

\bibitem{morDueG17}
Duersch, J.A., Gu, M.: Randomized {QR} with column pivoting.
\newblock {SIAM} J. Sci. Comput. \textbf{39}(4), C263--C291 (2017).
\newblock \doi{10.1137/15M1044680}

\bibitem{morEftKP11}
Eftang, J.L., Knezevic, D.J., Patera, A.T.: An hp certified reduced basis
  method for parametrized parabolic partial differential equations.
\newblock Math. Comput. Model. Dyn. Syst. \textbf{17}(4), 395--422 (2011).
\newblock \doi{10.1080/13873954.2011.547670}

\bibitem{Eve95}
Everson, R., Sirovich, L.: Karhunen--{L}o\`{e}ve procedure for gappy data.
\newblock J. Opt. Soc. Am. A \textbf{12}(8), 1657--1664 (1995)

\bibitem{gmsh09}
Geuzaine, C., Remacle, J.F.: Gmsh: A 3-{D} finite element mesh generator with
  built-in pre- and post-processing facilities.
\newblock Internat. J. Numer. Methods Engrg. \textbf{79}(11), 1309--1331
  (2009).
\newblock \doi{10.1002/nme.2579}

\bibitem{morGreMetal07}
Grepl, M., Maday, Y., Nguyen, N.C., Patera, A.T.: Efficient reduced-basis
  treatment of nonaffine and nonlinear partial differential equations.
\newblock {ESAIM}: Math. Model. Numer. Anal. \textbf{41}(3), 575--605 (2007).
\newblock \doi{10.1051/m2an:2007031}

\bibitem{morGre12}
Grepl, M.A.: Certified reduced basis methods for nonaffine linear time-varying
  and nonlinear parabolic partial differential equations.
\newblock Math. Models Methods Appl. Sci. \textbf{22}(3), 1150015, 40 (2012).
\newblock \doi{10.1142/S0218202511500151}

\bibitem{morGreP05}
Grepl, M.A., Patera, A.T.: A posteriori error bounds for reduced-basis
  approximations of parametrized parabolic partial differential equations.
\newblock {ESAIM}: Math. Model. Numer. Anal. \textbf{39}(1), 157--181 (2005).
\newblock \doi{10.1051/m2an:2005006}

\bibitem{morHaa17}
Haasdonk, B.: Reduced basis methods for parametrized {PDE}s---a tutorial
  introduction for stationary and instationary problems.
\newblock In: Model {R}eduction and {A}pproximation, \emph{Comput. Sci. Eng.},
  vol.~15, pp. 65--136. SIAM, Philadelphia, PA (2017).
\newblock \doi{10.1137/1.9781611974829.ch2}

\bibitem{morHaaDO11}
Haasdonk, B., Dihlmann, M., Ohlberger, M.: A training set and multiple bases
  generation approach for parameterized model reduction based on adaptive grids
  in parameter space.
\newblock Math. Comput. Model. Dyn. Syst. \textbf{17}(4), 423--442 (2011).
\newblock \doi{10.1080/13873954.2011.547674}

\bibitem{morHaaO08}
Haasdonk, B., Ohlberger, M.: Reduced basis method for finite volume
  approximations of parametrized linear evolution equations.
\newblock M2AN Math. Model. Numer. Anal. \textbf{42}(2), 277--302 (2008).
\newblock \doi{10.1051/m2an:2008001}

\bibitem{morHesRS16}
Hesthaven, J.S., Rozza, G., Stamm, B.: Certified Reduced Basis Methods for
  Parametrized Partial Differential Equations.
\newblock SpringerBriefs in Mathematics. Springer International Publishing
  (2016).
\newblock \doi{10.1007/978-3-319-22470-1}

\bibitem{morHesSZ14}
Hesthaven, J.S., Stamm, B., Zhang, S.: Efficient greedy algorithms for
  high-dimensional parameter spaces with applications to empirical
  interpolation and reduced basis methods.
\newblock {ESAIM}: Math. Model. Numer. Anal. \textbf{48}(1), 259--283 (2014).
\newblock \doi{10.1051/m2an/2013100}

\bibitem{morHesZ16}
Hesthaven, J.S., Zhang, S.: On the use of {ANOVA} expansions in reduced basis
  methods for parametric partial differential equations.
\newblock J. Sci. Comput. \textbf{69}(1), 292--313 (2016).
\newblock \doi{10.1007/s10915-016-0194-9}

\bibitem{morJiaC20}
Jiang, J., Chen, Y.: Adaptive greedy algorithms based on parameter-domain
  decomposition and reconstruction for the reduced basis method.
\newblock Internat. J. Numer. Methods Engrg. \textbf{121}(23), 5426--5445
  (2020).
\newblock \doi{10.1002/nme.6544}

\bibitem{morJiaCN17}
Jiang, J., Chen, Y., Narayan, A.: Offline-enhanced reduced basis method through
  adaptive construction of the surrogate training set.
\newblock J. Sci. Comput. \textbf{73}(2-3), 853--875 (2017).
\newblock \doi{10.1007/s10915-017-0551-3}

\bibitem{morMadS13}
Maday, Y., Stamm, B.: Locally adaptive greedy approximations for anisotropic
  parameter reduced basis spaces.
\newblock {SIAM} J. Sci. Comput. \textbf{35}(6), A2417--A2441 (2013).
\newblock \doi{10.1137/120873868}

\bibitem{Mah12}
Mahoney, M.W.: Algorithmic and statistical perspectives on large-scale data
  analysis.
\newblock In: Combinatorial Scientific Computing, Chapman \& Hall/CRC Comput.
  Sci. Ser., pp. 427--469. CRC Press, Boca Raton, FL (2012).
\newblock \doi{10.1201/b11644-17}

\bibitem{morMliGB15}
Mlinari{\'c}, P., Grundel, S., Benner, P.: Efficient model order reduction for
  multi-agent systems using {QR} decomposition-based clustering.
\newblock In: 54th IEEE Conference on Decision and Control (CDC), Osaka, Japan,
  pp. 4794--4799 (2015).
\newblock \doi{10.1109/CDC.2015.7402967}

\bibitem{morNar20}
Narayan, A.: Reduced order modeling and numerical linear algebra (2020).
\newblock
  \urlprefix\url{https://icerm.brown.edu/materials/Slides/sp-s20/Tutorial-Introductory_Talk_-_Snapshot-based_model_reduction_and_numerical_linear_algebra_part_1_]_Akil_Narayan,_University_of_Utah.pdf}.
\newblock ICERM Special Semester on Model and dimension reduction in uncertain
  and dynamic systems, Brown University, Providence, USA

\bibitem{morNegMA15}
Negri, F., Manzoni, A., Amsallem, D.: Efficient model reduction of parametrized
  systems by matrix discrete empirical interpolation.
\newblock J. Comput. Phys. \textbf{303}, 431--454 (2015).
\newblock \doi{10.1016/j.jcp.2015.09.046}

\bibitem{morNguRP09}
Nguyen, N.C., Rozza, G., Patera, A.T.: Reduced basis approximation and a
  posteriori error estimation for the time-dependent viscous {B}urgers'
  equation.
\newblock Calcolo \textbf{46}(3), 157--185 (2009).
\newblock \doi{10.1007/s10092-009-0005-x}

\bibitem{morTaineA15}
Paul-Dubois-Taine, A., Amsallem, D.: An adaptive and efficient greedy procedure
  for the optimal training of parametric reduced-order models.
\newblock Internat. J. Numer. Methods Engrg. \textbf{102}(5), 1262--1292
  (2015).
\newblock \doi{10.1002/nme.4759}

\bibitem{morPetal14}
Peherstorfer, B., Butnaru, D., Willcox, K., Bungartz, H.J.: Localized discrete
  empirical interpolation method.
\newblock {SIAM} J. Sci. Comput. \textbf{36}(1), A168--A192 (2014).
\newblock \doi{10.1137/130924408}

\bibitem{morPehDG18}
Peherstorfer, B., Drma\v{c}, Z., Gugercin, S.: Stability of discrete empirical
  interpolation and gappy proper orthogonal decomposition with randomized and
  deterministic sampling points.
\newblock e-prints 1808.10473v3, arXiv (2018).
\newblock \urlprefix\url{https://arxiv.org/abs/1808.10473v3}

\bibitem{morPehW15}
Peherstorfer, B., Willcox, K.: Online adaptive model reduction for nonlinear
  systems via low-rank updates.
\newblock {SIAM} J. Sci. Comput. \textbf{37}(4), A2123--A2150 (2015).
\newblock \doi{10.1137/140989169}

\bibitem{morPehW16}
Peherstorfer, B., Willcox, K.: Data-driven operator inference for nonintrusive
  projection-based model reduction.
\newblock Comp. Meth. Appl. Mech. Eng. \textbf{306}, 196--215 (2016).
\newblock \doi{10.1016/j.cma.2016.03.025}

\bibitem{morPehZB13}
Peherstorfer, B., Zimmer, S., Bungartz, H.J.: Model reduction with the reduced
  basis method and sparse grids.
\newblock In: Sparse grids and applications, \emph{Lect. Notes Comput. Sci.
  Eng.}, vol.~88, pp. 223--242. Springer, Heidelberg (2013).
\newblock \doi{10.1007/978-3-642-31703-3}

\bibitem{morQuaMN16}
Quarteroni, A., Manzoni, A., Negri, F.: {R}educed {B}asis {M}ethods for
  {P}artial {D}ifferential {E}quations, \emph{La Matematica per il 3+2},
  vol.~92.
\newblock Springer International Publishing (2016).
\newblock \doi{10.1007/978-3-319-15431-2}

\bibitem{morQuaRM11}
Quarteroni, A., Rozza, G., Manzoni, A.: Certified reduced basis approximation
  for parametrized partial differential equations and applications.
\newblock J. Math. Ind. \textbf{1}, Art. 3, 44 (2011).
\newblock \doi{10.1186/2190-5983-1-3}

\bibitem{morRavS20}
Rave, S., Saak, J.: A non-stationary thermal-block benchmark model for
  parametric model order reduction.
\newblock e-print 2003.00846, arXiv (2020).
\newblock \urlprefix\url{https://arxiv.org/abs/2003.00846}.
\newblock Math.NA

\bibitem{morRozHS20}
Rozza, G., Hess, M., Stabile, G., Tezzele, M., Ballarin, F.: Basic ideas and
  tools for projection-based model reduction of parametric partial differential
  equations.
\newblock In: P.~Benner, S.~Grivet-Talocia, A.~Quarteroni, G.~Rozza,
  W.~Schilders, L.M. Silveira (eds.) Snapshot-{B}ased {M}ethods and
  {A}lgorithms, pp. 1--47. De Gruyter (2021).
\newblock \doi{10.1515/9783110671490-001}

\bibitem{morRozHP08}
Rozza, G., Huynh, D.B.P., Patera, A.T.: Reduced basis approximation and a
  posteriori error estimation for affinely parametrized elliptic coercive
  partial differential equations: application to transport and continuum
  mechanics.
\newblock Arch. Comput. Methods Eng. \textbf{15}(3), 229--275 (2008).
\newblock \doi{10.1007/s11831-008-9019-9}

\bibitem{morSai20}
Saibaba, A.K.: Randomized discrete empirical interpolation method for nonlinear
  model reduction.
\newblock {SIAM} J. Sci. Comput. \textbf{42}(3), A1582--A1608 (2020).
\newblock \doi{10.1137/19M1243270}

\bibitem{morSen08}
Sen, S.: Reduced-basis approximation and a posteriori error estimation for
  many-parameter heat conduction problems.
\newblock Numerical Heat Transfer, Part B: Fundamentals \textbf{54}(5),
  369--389 (2008).
\newblock \doi{10.1080/10407790802424204}

\bibitem{morTezBR18}
Tezzele, M., Ballarin, F., Rozza, G.: Combined parameter and model reduction of
  cardiovascular problems by means of active subspaces and {POD}-{G}alerkin
  methods.
\newblock In: Mathematical and numerical modeling of the cardiovascular system
  and applications, \emph{SEMA SIMAI Springer Ser.}, vol.~16, pp. 185--207.
  Springer, Cham (2018).
\newblock \doi{10.1007/978-3-319-96649-6_8}

\bibitem{Zhaetal}
Zha, H., He, X., Ding, C., Gu, M., Simon, H.D.: Spectral relaxation for k-means
  clustering.
\newblock In: T.G. Dietterich, S.~Becker, Z.~Ghahramani (eds.) Advances in
  Neural Information Processing Systems 14, pp. 1057--1064. MIT Press (2002).
\newblock
  \urlprefix\url{http://papers.nips.cc/paper/1992-spectral-relaxation-for-k-means-clustering.pdf}

\bibitem{morZhaetal15}
Zhang, Y., Feng, L., Li, S., Benner, P.: An efficient output error estimation
  for model order reduction of parametrized evolution equations.
\newblock {SIAM} J. Sci. Comput. \textbf{37}(6), B910--B936 (2015).
\newblock \doi{10.1137/140998603}

\end{thebibliography}

\end{document}